\documentclass[review,10pt]{elsarticle}
\usepackage[utf8]{inputenc}
\usepackage{microtype}
\usepackage[margin=2.5cm]{geometry}
\usepackage{setspace}
\usepackage[none]{hyphenat}
\usepackage[pagewise]{lineno}
\usepackage[english]{babel}
\usepackage{color}
\usepackage{booktabs}
\usepackage{longtable}
\usepackage[most]{tcolorbox}
\usepackage{ragged2e}
\usepackage[justification=justified]{caption}
\usepackage{enumitem}

\usepackage{amsmath,amsfonts,amssymb,mathtools}
\allowdisplaybreaks 
\usepackage{cases}
\usepackage{multirow}
\usepackage{graphicx,float}
\usepackage{algpseudocode}
\usepackage{subcaption}
\usepackage{eurosym}
\usepackage{amstext} 
\usepackage{arydshln}
\usepackage{flushend}
\usepackage{threeparttable}
\usepackage{url}
\usepackage{algorithm}

\usepackage[]{hyperref}
\hypersetup{colorlinks, linkcolor=blue, citecolor=blue, urlcolor=blue}

\DeclareRobustCommand{\officialeuro}{%
  \ifmmode\expandafter\text\fi
  {\fontencoding{U}\fontfamily{eurosym}\selectfont e}}
\usepackage{mdframed}
\immediate\write18{%
makeindex Main_File.nlo -s nomencl.ist -o Main_File.nls 
}
\usepackage{multicol} 
\usepackage{nomencl}[noprefix] 
\makenomenclature
\renewcommand*\nompreamble{\begin{multicols}{2}}
\renewcommand*\nompostamble{\end{multicols}}

\def\tsc#1{\csdef{#1}{\textsc{\lowercase{#1}}\xspace}}
\tsc{WGM}
\tsc{QE}
\tsc{EP}
\tsc{PMS}
\tsc{BEC}
\tsc{DE}

\begin{document}

    \justifying
    \begin{frontmatter}   
        \title{Hybrid Adaptive Robust-Stochastic Optimization Model for the Design of a Photovoltaic–Battery Energy Storage System }
    
        \author[1,3]{Alba Lun Mora Pous\corref{mycorrespondingauthor}}
        \cortext[mycorrespondingauthor]{Corresponding author}
        \ead{alba.lun.mora@estudiantat.upc.edu}
        \author[1,2]{Fernando Garcia-Muñoz}
        \author[1,2]{Natalia Jorquera-Bravo}
        \author[1]{Ricardo Aranguiz}
        \author[4]{Valentina Bugueño Olivos}

        \address[1]{University of Santiago of Chile (USACH), Faculty of Engineering, Industrial Engineering Department, Chile}
        \address[2]{University of Santiago of Chile (USACH), Faculty of Engineering, Program for the Development of Sustainable Production Systems (PDSPS), Chile}
        \address[3]{Universitat Politècnica de Catalunya (UPC), Escola Tècnica Superior d'Engyinyeria Industrial de Barcelona, Spain}
        \address[4]{University of Santiago of Chile (USACH), Faculty of Engineering, Electrical Engineering Department, Chile}

    \begin{abstract}
       Future energy projections and their inherent uncertainty play a key role in the design of photovoltaic–battery energy storage systems (PV-BESS) for household use. In this study, both stochastic and robust optimization techniques are simultaneously integrated into a Hybrid Adaptive Robust–Stochastic Optimization (HARSO) model. Uncertainty in future PV generation is addressed using a stochastic approach, while uncertainty in power demand is handled through robust optimization. To solve the tri-level structure emerging from the hybrid approach, a Column-and-Constraint Generation (CCG) algorithm is implemented. 
       The model also accounts for battery degradation by considering multiple commercially available battery chemistries, enabling a more realistic evaluation of long-term system costs and performance. To demonstrate its applicability, the model is applied to a case study involving the optimal design of a PV-BESS system for a household in Spain. The empirical analysis includes both first-life (FL) and second-life (SL) batteries with different chemistries, providing a comprehensive evaluation of design alternatives under uncertainty. Results indicate that the optimal solution is highly dependent on the level of robustness considered, leading to a shift in design strategy. Under less conservative settings, robustness is achieved by increasing battery capacity, while higher levels of conservatism favor expanding PV capacity to meet demand. 
       
    \end{abstract}


\end{frontmatter}

\section{Introduction}\label{Intro_Chap}

CO\textsubscript{2} emissions have increased by more than 100\% over the past 50 years, with energy-related emissions rising by approximately 900 Mt between 2019 and 2023~\cite{58}, posing not only environmental and planetary threats but also significant risks to human health~\cite{59}. Growing environmental and social concerns have accelerated the shift toward cleaner energy solutions, such as photovoltaic (PV) panels for household power. Adopting this technology can benefit users in many ways. In addition to environmental benefits, households become more autonomous from the electric grid, and electricity bills can be reduced~\cite{PVBenefits}. PV systems are often paired with Battery Energy Storage Systems (BESS), forming a combined PV-BESS setup. This combined setup is more beneficial than standalone PV installations, they can increase self-consumption by storing surplus energy, provide backup power during outages, and further reduce electricity bills by purchasing or selling energy when prices are most cost-effective \cite{PV-BESSBenefits}. Moreover, PV-BESS configurations open the door to reusing second-life (SL) batteries, retired electric vehicle (EV) batteries no longer suitable for automotive use, which can be repurposed for stationary applications. This not only extends battery life cycles but also reduces waste and production-related emissions~\cite{18}.

However, to fully realize the benefits of a PV-BESS system, proper sizing is essential to ensure cost-effectiveness and justify the initial investment. Optimal sizing requires to understand its operating strategy and how the BESS will perform over time. To do this, several aspects must be considered, such as electricity purchase and selling prices, the capacity of the PV panels, the performance of the battery modules, and expected patterns of energy generation and use. Several studies have addressed optimal sizing strategies for such systems; a summary of relevant work is provided in ~\cite{85}.

A major challenge in this process is dealing with the uncertainties that affect energy systems~\cite{108}. On the production side, PV output depends heavily on local weather conditions, particularly solar irradiance, which is difficult to predict accurately over long time horizons. While no energy is generated during nighttime hours, daytime production remains highly variable due to factors such as cloud cover 
and seasonal variations. This leads to substantial differences in solar energy availability throughout the year, adding another layer of complexity to production forecasting \cite{110}. On the consumption side, 
household energy demand is even harder to estimate. Residential energy use can fluctuate significantly over time due to changes in user behavior, occupancy patterns, the introduction of new appliances, regulations, or even lifestyle changes~\cite{uncertaintyDemand}. 
These uncertainties must be properly accounted for to ensure a reliable and robust system design. Another critical aspect of PV-BESS optimization is selecting the most suitable battery type. Battery technologies vary in terms of performance and degradation behavior~\cite{63}, and SL batteries typically degrade faster than first-life (FL) alternatives~\cite{6,12,15,16}. Incorporating degradation dynamics into the modeling framework is thus vital to improving model realism and predictive accuracy.

 Given these considerations, it becomes relevant to develop PV-BESS sizing models that not only incorporate battery degradation over time, but also account for energy uncertainties. Such models are essential to obtain reliable and realistic system configurations that remain robust under evolving conditions. To address this challenge, this study proposes a modeling approach that explicitly considers both energy demand and PV generation uncertainty, battery degradation and different performance depending on battery chemistry.

The remainder of the article is organized as follows. Section~\ref{sec:literature} reviews existing approaches for PV-BESS sizing. Section~\ref{sec:methodology} outlines the methodology adopted in this work to estimate electricity pricing, future consumption patterns, uncertainty bounds, scenario generation, and battery degradation. Section~\ref{sec:Opt_model} introduces the optimization model and the proposed approaches to handle uncertainty. Section \ref{sec:Algorithm} introduces the solution strategy for addressing the tri-level structure. Section~\ref{sec:case study} presents a case study and discusses the computational results. Finally, Section~\ref{sec:conclusions} summarizes the main findings and conclusions of the study.

\section{Literature review}\label{sec:literature}

The following literature review aims to synthesize previous research on PV-BESS sizing approaches, highlighting the strategies used to characterize uncertainty.

  Deterministic approaches have been widely used to determine the optimal size of PV-BESS systems. These methods neglect the inherent uncertainty associated with BESS and the PV generation, using daily averages of solar energy and load demand, making the assumption that they are sufficiently representative over time \cite{88}. In fact, as highlighted in \cite{102}, the number of studies adopting a deterministic framework significantly outweighs those that explicitly address uncertainty. In
\cite{65}, a Mixed Integer Linear Programming (MILP) model is proposed to optimize PV-BESS sizing that takes into consideration important factors such as optimal PV tilt angle, vehicle-to-home (V2H) availability, battery degradation, and even load scheduling of controllable equipment. The authors of \cite{66} employ a Genetic Algorithm (GA) to develop an online tool for determining the optimal sizing of PV-BESS systems based on the household's self-sufficiency expectations. The approach relies on a database of precomputed historical simulations, assuming these are sufficiently representative to guide new cases. This enables the tool to deliver real-time recommendations without solving a new optimization problem each time. In \cite{73}, the authors solve the optimization problem considering SL batteries, with emphasis on the degradation that they might suffer over the years. Other studies also employ deterministic approaches to optimally design a PV-BESS system for communities involving several households, such as in \cite{63-a}. While deterministic optimization offers advantages such as lower computational burden and simpler formulations, its underlying assumptions can lead to suboptimal or unreliable solutions, especially in the medium and long term, where uncertainty plays a critical role.

There exist different approaches to uncertainty assessment in energy system optimization models, and the two main approaches are Stochastic Optimization (SO) \cite{birge2011introduction} and Robust Optimization (RO) \cite{ben2009robust}. Stochastic problems model uncertainty  by discretizing known probability distributions into scenarios. In RO, however, the uncertain parameters have set-based definitions and require minimal uncertainty information, which has the advantage that only a range of variation is required for each parameter and no probability distribution is needed.

Two-stage stochastic optimization has been applied in energy system planning problems, as its structure naturally aligns with the sequential nature of investment and operation decisions. In this framework, long-term decisions such as the sizing of storage and conversion technologies are made in the first stage, while operational strategies are optimized in the second stage under a set of predefined scenarios. This approach is adopted in \cite{200}, where the optimal sizing of a BESS is determined based on four scenarios representing annual wind generation profiles. However, the focus remains primarily on generation variability, with limited consideration of consumption-side uncertainty. Similar formulations are used in \cite{9408593,LI20211837,LIN2024112862}, where scenario-based investment and operation models are developed for different configurations. For example, \cite{9408593} considers uncertainty in electricity demand, heat demand, and market prices, while \cite{LI20211837} jointly optimizes the sizing and dispatch of a hybrid wind–concentrated solar power system with electric heaters. In turn, \cite{LIN2024112862} evaluates multiple PV-BESS-hydrogen configurations, balancing economic, reliability, and sustainability objectives through a modified genetic algorithm coupled with a decision-making technique based on distance to an ideal solution. While the two-stage structure offers conceptual and computational advantages, a known limitation of this approach is the difficulty in accurately characterizing the probability distribution functions of long-term uncertainties, which often leads to questionable representativeness of the scenario set.

  The limitation of the above two-stage stochastic approach motivates using alternative approaches, such as RO, which rely on uncertainty sets rather than probabilistic assumptions. Thus, RO comes in several model types, depending on how uncertainty is represented, how conservative the solution is, and what trade-offs are allowed. In \cite{99}, they apply an RO to household load scheduling, aiming to determine the optimal operation strategy of the BESS. This study considers uncertainty exclusively in electricity prices and does not address the sizing of the PV-BESS system. In contrast, \cite{100} also focuses on optimizing the BESS operation rather than sizing, but it incorporates uncertainty in both energy demand and renewable energy production. These uncertainties are modeled under the same uncertainty framework by using interval uncertainty within a set of budget of uncertainty (BoU), known as Bertsimas \& Sim approach \cite{bertsimas2004price}. 


Other RO approaches have been made in \cite{CORREAFLOREZ2018433}, where uncertainty is considered in prices, energy load, and PV energy production by using an Adaptive Robust Optimization (ARO) approach. Although this approach is computationally more demanding, it enables solutions that adapt after the realization of uncertainty \cite{sun2021robust}. This is particularly advantageous, as it leads to less conservative and more realistic decisions, which is something that classical RO does not allow. Thus, authors from \cite{95} focus on optimizing the system operation for a community and its interaction with electricity markets. An ARO model is presented, and a comparative analysis of different modifications to this formulation is conducted to assess its potential and limitations. In the initial approach, uncertainty is captured through a single variable: the net load, defined as the difference between consumption and PV generation. This results in a more compact formulation that may appear conservative. However, when comparing with alternative models that treat generation and consumption uncertainties separately, it is observed that the latter approach offers greater flexibility in tuning robustness levels, which can reduce conservatism and improve performance in certain cases. Alternatively, in \cite{109}, an ARO model for optimal design is presented, where uncertainty in both electricity generation and demand is modeled using user-defined bounded intervals within polyhedral uncertainty sets. 

Other studies address the optimization problem using a hybrid model that integrates both RO and SO to handle uncertainties. Such hybrid approaches can outperform standalone RO or SO by balancing overly conservative sizing and overly optimistic planning \cite{sun2021robust}. For instance, \cite{106} adopts this strategy by incorporating uncertainties in both energy prices and PV generation. For day-ahead energy prices, a worst-case scenario is contemplated, while probable scenarios are used to model the uncertainties of PV generation and real-time market prices. However, the study focuses solely on optimal energy management, without addressing system sizing. Moreover, it assumes a deterministic household load profile and does not account for battery degradation over time. Similarly, \cite{107} also focuses on finding the optimal energy management strategy for a household equipped with multiple types of loads and energy sources. In this case, a robust method to tackle uncertainty is used to estimate PV generation and real time energy prices. On the other hand, stochastic optimization is applied to the charging and discharging behavior of EVs, as the household is assumed to operate under a vehicle-to-grid (V2G) framework. However, since the scope differs from that of optimal PV-BESS sizing, the study does not address BESS degradation. Additionally, the household demand is modeled as the aggregate of individual load components, each represented separately and considered deterministic. 

As highlighted in the referenced studies, employing different modeling approaches within a unified optimization framework is essential, since the nature of uncertainty varies widely depending on several factors. On the one hand, it depends on the tolerance of the decision-maker. Typically, SO seeks to optimize the expected value, while RO aims at safeguarding against the worst case. In this sense, the choice between SO and RO reflects a strategic preference in handling uncertainty.
On the other hand, it is also determined by the characteristics of the modeled parameter. For instance, PV generation lends itself to relatively accurate probabilistic modeling, as it follows well-understood physical patterns, no generation occurs at night, and peak output aligns with solar position due to Earth's rotation and axial tilt~\cite{earthRotation}. In contrast, electrical demand is driven by human behavior, which is inherently volatile and difficult to predict. Although it can also be modeled probabilistically, its characterization through probability distributions is considerably more challenging and generally less reliable. This challenge is expected to intensify in the coming years due to electrification trends, such as the adoption of electric vehicles and electric heating, which will significantly alter demand profiles and increase bot temporal and spatial variability~\cite{understandingUncertainty}. 

  
This paper aims to address the gap in the literature on the optimal design of PV-BESS systems using a Hybrid Adaptive Robust-Stochastic Optimization (HARSO) approach, which models uncertainty in PV generation using SO and in energy demand using RO. While several studies have explored PV-BESS sizing, they typically adopt either a deterministic, robust, or stochastic formulation, with limited attention to hybrid models. Additionally, the BESS is explicitly modeled to account for degradation over time and its impact on system performance. The model includes the selection among different battery technologies, such as SL and FL options, allowing the optimization to endogenously determine the most suitable technology based on trade-offs between cost and degradation behavior.


\section{Methodology}\label{sec:methodology}
\begin{figure}[htbp] 
    \centering
    \includegraphics[width=1.0\textwidth]{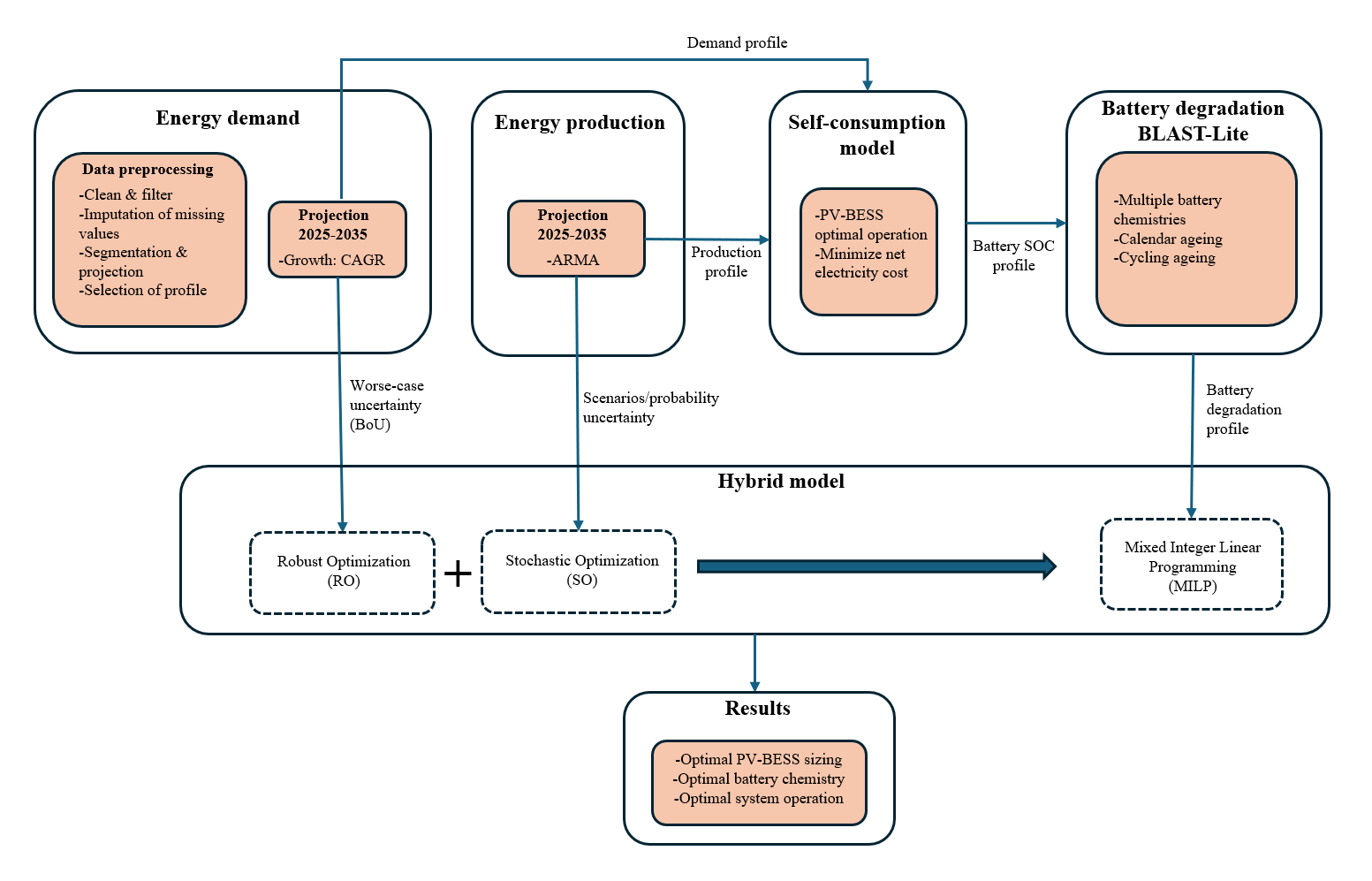} 
    \caption{Methodology map.}
    \label{fig:methodology_map}
\end{figure}

This section describes the methodological steps to build a consistent and data-driven representation of the PV-BESS system's future behaviour, encompassing energy demand, PV generation, and battery degradation. This will be used as inputs parameters for the HARSO model presented in Section \ref{sec:Opt_model}.
Figure~\ref{fig:methodology_map} presents a visual summary of the process, and each step is detailed below.


The first task focuses on forecasting electricity demand with hourly resolution for the coming years. Historical household consumption data are used to generate normalized hourly demand profiles for each year, which define annual demand intervals and serve as the basis for the RO component through the BoU framework. To extend these projections, we apply the Compound Annual Growth Rate (CAGR) formula, a method widely used in energy forecasting \cite{24} because it provides consistent growth estimates even when the data exhibit variability \cite{23}. In parallel, PV generation is represented through a stochastic Auto-Regressive Moving Average (ARMA) process, selected for its ability to capture the dependence of solar output on past values \cite{JI2011808}. The model combines an autoregressive component, which relates current generation to previous observations, with a moving-average component, which accounts for the effect of past random fluctuations, complemented by a noise term to reflect residual uncertainty. A simple ARMA structure is adopted as a trade-off between tractability and the capacity to reproduce short-term variability. Seasonal and interannual differences are incorporated through a set of representative daily profiles constructed for each projected year.

Building on these inputs, a simplified deterministic optimization model is first employed to estimate the annual operating profile of a given BESS configuration, which is then used to evaluate battery degradation. The resulting operational pattern is processed in a MATLAB-based degradation model \cite{BLAST} to quantify yearly performance losses across different battery technologies, including both FL and SL options. With the projected demand intervals, PV generation scenarios, and technology-specific degradation curves as inputs, the HARSO model is finally solved to determine the optimal PV-BESS sizing and to assess its long-term performance under uncertainty.

\section{Optimization model and uncertainty}\label{sec:Opt_model}
In this section, we first describe a MILP model to address the optimal sizing and operational planning of a residential PV-BESS system within a deterministic framework. This model is then extended to incorporate uncertainty, resulting in the HARSO formulation. All parameters and variables are defined progressively as the corresponding equations are introduced, while a complete summary is also provided in Table~\ref{Nomenclature_table} for ease of reference.  All variables are continuous and non-negative, except those explicitly indicated as binary.

\subsection{Deterministic Problem}\label{Det_Model_Subsection}
The proposed optimization model addresses the optimal sizing and operational planning of a residential PV-BESS system over a multi-year planning horizon \( \mathcal{Y} \), with hourly time resolution \( t \in \mathcal{T} \). The model has to determine whether investment in a PV system is profitable and, if so, it has to select at most one battery technology from the set of candidates \( j \in \mathcal{J} \). Each battery is characterized by different investment costs and degradation profiles. The objective is to determine the optimal PV and BESS capacities that minimize the total cost, which includes both investment and operational expenses associated with PV generation and battery cycling. To avoid overinvestment, the system operates in grid-connected mode, allowing the user to exchange electricity with the main grid at each time step. In the case of a self-consumption deficit, electricity can be purchased from the grid as \( p^{bg}_{t,y} \) at a price \( \lambda^{bg}_{t} \), and in periods of energy surplus, electricity can be sold back to the grid as \( p^{sg}_{t,y} \) at a price \( \lambda^{sg}_{t} \). The complete optimization problem is formulated as a Mixed-Integer Linear Programming (MILP) model, with binary variables governing investment and operational mode decisions. Thus, the objective function is defined as follows: 
\begin{align}\label{FO_Det_model}
\min\quad z =\; & C^{PV}\, \gamma^{pv}
  \;+\;\sum_{j\in \mathcal{J}}C_{j}^{BT}\,\gamma^{bt}_{j}
  \;+\;\sum_{t\in \mathcal{T}}\sum_{y\in \mathcal{Y}}\Bigl(\lambda^{bg}_{t} \,p^{bg}_{t,y}
    -\lambda^{sg}_{t} \,p^{sg}_{t,y}
    +pg_{t,y} \,C_{\mathrm{op}}^{PV}\Bigr)\notag\\
  &\quad+\;\sum_{j\in \mathcal{J}}\sum_{t\in \mathcal{T}}\sum_{y\in \mathcal{Y}}\,ds_{j,t,y}  \,C_{\mathrm{op}}^{BT}
\end{align}

The objective function in~\eqref{FO_Det_model} minimizes the total system cost over the entire planning horizon. It includes both investment and operational components. The first two terms represent the investment costs, such that \( C^{PV} \) is the marginal cost per kW of installed PV capacity \( \gamma^{pv} \), and \( C^{BT}_{j} \) is the cost per kWh for each candidate battery technology \( j \in \mathcal{J} \), associated with its installed capacity \( \gamma^{bt}_{j} \). The third term represents the cumulative operational cost over all periods \( t \in \mathcal{T} \) and years \( y \in \mathcal{Y} \). It includes the cost of electricity imports from the grid \( p^{bg}_{t,y} \) at price \( \lambda^{bg}_{t} \), the revenue from electricity exports to the grid \( p^{sg}_{t,y} \) at price \( \lambda^{sg}_{t} \), and a linear operational cost \( C_{\mathrm{op}}^{PV} \) associated with PV generation \( pg_{t,y} \), reflecting the economic uses of the PV asset. The last term captures the battery-related operational cost, modeled as the cost \( C_{\mathrm{op}}^{BT} \) per unit of discharged energy \( ds_{j,t,y} \) associated to each of the battery technology \( j \in \mathcal{J} \). 

  The feasible solution space is constrained by the following set of operational and investment constraints, which govern power balance, PV generation limits, BESS dynamics, and technology selection decisions:
\begin{subequations}\label{Det_Const_model}
\begin{align}
&pg_{t,y}\;-\;PL_{t,y}\;-\;p^{sg}_{t,y}\;+\;p^{bg}_{t,y}
   \;+\;\sum_{j\in \mathcal{J}}\bigl(ds_{j,t,y}-ch_{j,t,y}\bigr)
   = 0
   && \forall\,t\in \mathcal{T},\,y\in \mathcal{Y}
   \label{Balance_Const}\\
&pg_{t,y}\;\le\;\gamma^{pv}\;\cdot\;\overline{PG}_{t,y}
   && \forall\,t\in \mathcal{T}, y\in \mathcal{Y} \label{PV_generation_Const}\\
&\gamma^{pv}\;\le\;\overline{\Gamma^{PV}} \label{PV_Cap_Const}\\
&\begin{aligned}
\text{soc}_{j,t,y} =
\begin{cases}
SOC^{i} \cdot \gamma^{bt}_{j} + \varphi \cdot ch_{j,t,y} - \dfrac{1}{\varphi} \cdot ds_{j,t,y} & \text{if } t = 1 \\
SOC^{l} \cdot \gamma^{bt}_{j} & \text{if } t = 24 \\
\text{soc}_{j,t-1,y} + \varphi \cdot ch_{j,t,y} -  \dfrac{1}{\varphi} \cdot ds_{j,t,y} & \text{otherwise}
\end{cases}
\end{aligned}
  && \forall\,j\in \mathcal{J},t\in \mathcal{T},\,\,y\in \mathcal{Y} \label{SOC_Const}\\
&DG_{j,y}\cdot\gamma^{bt}_{j}\cdot\,\underline{SOC}
   \;\le\;soc_{j,t,y}
   \;\le\;DG_{j,y}\cdot\,\gamma^{bt}_{j}\cdot\,\overline{SOC}
   && \forall\,j\in \mathcal{J},t\in \mathcal{T},\,y\in \mathcal{Y} \label{SOC_Limit_Const}\\
&ch_{j,t,y}\;\le\;PB\;\cdot\;w_{j,t,y}
   && \forall\,j\in \mathcal{J},t\in \mathcal{T},\,y\in \mathcal{Y}\label{CH_Const}\\
&ds_{j,t,y}\;\le\;PB\,(\nu^{bt}_{j} - w_{j,t,y})
   && \forall\,j\in \mathcal{J},t\in \mathcal{T},y\in \mathcal{Y}\label{DS_Const}\\
&w_{j,t,y}\;\le\;\nu^{bt}_{j}
   && \forall\,j\in \mathcal{J},t\in \mathcal{T},y\in \mathcal{Y}\label{Binary_w_Const}\\
&\gamma^{bt}_{j}\;\le\;\nu^{bt}_{j}\;\cdot\;\overline{\Gamma^{BT}}
   && \forall\,j\in \mathcal{J\label{BESS_Cap_Const}} \displaybreak[1]\\[4pt]
    &\sum_{j\in \mathcal{J}}\nu^{bt}_{j}\;\le\;1 && \label{BESS_Tech_Const}
\end{align}
\end{subequations}

  The power balance constraint~\eqref{Balance_Const} ensures that, at every time step \( t \in \mathcal{T} \) and year \( y \in \mathcal{Y} \), the total energy demand \( PL_{t,y} \) is met by PV generation \( pg_{t,y} \), grid imports \( p^{bg}_{t,y} \), and BESS discharging \( ds_{j,t,y} \), while accounting for grid exports \( p^{sg}_{t,y} \) and BESS charging \( ch_{j,t,y} \). The PV generation constraint~\eqref{PV_generation_Const} links the generated power \( pg_{t,y} \) to the installed capacity \( \gamma^{pv} \) and the irradiance profile \( \overline{PG}_{t,y} \), while the investment constraint~\eqref{PV_Cap_Const} enforces an upper bound \( \overline{\Gamma^{PV}} \) on the installed PV capacity. Likewise, BESS operation is governed by the SOC update equation~\eqref{SOC_Const}, which accounts for charging and discharging efficiencies \( \varphi \) and includes boundary conditions for the first and last hour of each day. The SOC is further bounded in~\eqref{SOC_Limit_Const} between minimum and maximum allowable levels, scaled by the degradation factor \( DG_{j,y} \) and the battery capacity \( \gamma^{bt}_j \). The binary variable \( w_{j,t,y} \) determines whether the battery \( j \) is in charging mode, thus, Constraints~\eqref{CH_Const} and~\eqref{DS_Const} enforce mutual exclusivity between charging and discharging and prevent operation if the battery has not been selected (\( \nu^{bt}_{j} = 0 \)). Constraint~\eqref{Binary_w_Const} ensures that the operational mode is only active if the corresponding battery is installed. Finally, the installation constraint~\eqref{BESS_Cap_Const} sets an upper bound on the installed capacity of each battery technology \( j \in \mathcal{J} \), conditional on the binary selection variable \( \nu^{bt}_j \), and constraint~\eqref{BESS_Tech_Const} restricts the solution to select at most one battery technology for installation.

\subsection{Handling uncertainty}\label{Section_HARSO_model}
  The deterministic optimization model presented above assumes perfect knowledge of key parameters such as electricity demand, PV generation, electricity prices, and technology costs over the entire planning horizon. While this assumption simplifies the formulation, it fails to capture the inherent variability and unpredictability of real-world energy systems. In practice, assuming full certainty in such long-term investment planning problems leads to solutions that may be either infeasible or suboptimal under actual operational conditions. 
  
  In this work, we address uncertainty exclusively in electricity demand and PV generation, while keeping electricity prices and technology costs as known parameters. The inclusion of uncertainty in price parameters is considered beyond the scope of this study and is left for future extensions. It is important to note that the nature of uncertainty differs significantly between PV generation and electricity demand. Although PV output can exhibit short-term variability due to weather conditions, its long-term behavior typically follows well-defined seasonal and diurnal patterns, which allows for a probabilistic characterization based on historical data. In contrast, residential electricity demand is more difficult to predict over long horizons. It is shaped by complex and evolving factors such as lifestyle changes, appliance usage, and the ongoing electrification associated with the energy transition. These factors introduce a higher degree of uncertainty that is harder to capture using a single probability distribution. 
  
  Given these differences, we adopt a hybrid uncertainty modeling framework that combines stochastic and robust optimization. PV generation uncertainty is modeled stochastically using a scenario-based approach derived from time-series models. Meanwhile, electricity demand uncertainty is represented using a robust formulation with an uncertainty set, which allows the model to remain conservative with respect to worst-case variations while maintaining tractability. This hybrid treatment enables a more realistic and structurally accurate representation of the system's uncertainties, aligning the optimization model with the distinct behaviors of its key components.

Following the structure presented in~\cite{sun2021robust}, the optimization problem used in this study can be represented as a tri-level HARSO model, formulated as:
\begin{equation}
\min_{x \in \mathcal{X}} \ \max_{u \in \mathcal{U}} \ \min_{y \in \mathcal{Y}(x,u,\xi)} \ \mathbb{E}_{\xi} \left[ f(x, u, y, \xi) \right].
\end{equation}

In this formulation, \( x \) denotes the planning (here-and-now) decisions, such as the sizing of PV and BESS assets; \( u \in \mathcal{U} \) denotes deviations in electricity demand as captured by the uncertainty set $\mathcal{U}$; \( \xi \) corresponds to a random variable with a known probability distribution for the PV generation; 
and \( y \in \mathcal{Y}(x, u, \xi) \) includes the recourse variables adjusted ex-post, such as battery dispatch and load supplied. This tri-level structure captures the hybrid nature of uncertainty by embedding both robust and stochastic components: the outer maximization over \( u \) accounts for the worst-case demand realization, while the expectation operator \( \mathbb{E}_{\xi} \) reflects the probabilistic behavior of PV generation. The adaptive recourse minimizes the expected operational cost for each scenario, conditional on both the investment plan and the worst-case realization of demand. Additionally, this work adopts a BoU to characterize the uncertainty set associated with long-term uncertainty in the robust component of the  HARSO model. Thus, following the formulation proposed in~\cite{ATTARHA2018132}, the uncertainty set is defined as:
\begin{equation}
\mathcal{U} = \left\{ u_t = \bar{u}_t + \hat{u}_t z_t : z_t \in [-1,1], \ \sum_{t \in \mathcal{T}} |z_t| \leq \Gamma \right\},
\end{equation}
where \( \bar{u}_t \) denotes the nominal value of the uncertain parameter at time \( t \), \( \hat{u}_t \) is the maximum allowed deviation, \( z_t \) are auxiliary variables bounded within \( [-1,1] \), and \( \Gamma \) is a scalar parameter that limits the total deviation over the planning horizon. This structure enables the model to account for uncertainty in a controlled manner by restricting the number and intensity of deviations from nominal values. When \( \Gamma = 0 \), the model reduces to the deterministic case, while \( \Gamma = |\mathcal{T}| \) corresponds to a fully conservative approach that protects against the worst-case realization across all periods.

On the other hand, the random variable $ \xi$ is approximated using a \textit{Sample Average Approximation} (SAA) approach \cite{SAA}. The SAA approach replace the original distribution with a sample of scenarios $\mathcal{S}$, where each scenario $s \in \mathcal{S}$ correspond to a realization of $\xi$ with positive probability $\rho_s>0$ such that $\sum_{s\in\mathcal{S}}\rho_s=1$.  

  Finally, connecting the compact HARSO structure with the deterministic formulation introduced in Section~\ref{Det_Model_Subsection}, the first-stage decision variables \( x \in \mathcal{X} \) correspond to the investment-related decisions, namely the sizing of PV and BESS assets \( \gamma^{pv} \) and \( \gamma^{bt}_j \), as well as the binary selection variables \( \nu^{\mathrm{bt}}_j \) associated with the installation of each battery type \( j \). The long-term uncertainty vector \( u \in \mathcal{U} \) is mapped to the hourly demand profiles \( PL_{y,t} \), whose deviations are modeled using the BoU described above. The short-term stochastic component \( \xi \) captures the variability of PV generation \( \overline{PG}_{t,y} \) across scenarios, and the second-stage recourse include operational decisions related to dispatch, battery control, and energy balancing $y =\{  pg_{t,y}, p^{sg}_{t,y}, p^{bg}_{t,y}, ds_{j,t,y}, ch_{j,t,y}, soc_{j,t,y}, w_{j,t,y}   \}$.




\section{Algorithm}\label{sec:Algorithm}
The HARSO model presented in Section~\ref{Section_HARSO_model} exhibits a nested tri-level structure, where investment decisions are followed by an adversarial worst-case realization of uncertain parameters, and subsequently by a scenario-based recourse stage \cite{sun2021robust}. Solving such models directly using commercial solvers is computationally intractable due to the non-convex structure induced by the min-max-min formulation. Specifically, the robust component requires identifying the worst-case realization within a high-dimensional uncertainty set, while the stochastic part involves computing expectations over multiple scenarios. This combination leads to a non-convex problem that cannot be reformulated into a tractable single-level equivalent. 

To address this complexity, decomposition techniques such as Benders decomposition and Column-and-Constraint Generation (CCG) have been proposed as effective methods for solving robust formulations~\cite{113}. The main difference between these approaches lies in the structure of the Master Problem (MP) and the nature of the cuts added at each iteration. In Benders decomposition, the MP includes only the first-stage variables and is updated with dual cuts derived from the subproblem (SP). In contrast, CCG constructs an extended MP that incorporates both first- and second-stage variables, but only the second-stage variables associated with the worst-case scenario identified in the current iteration. The cuts in CCG are primal and correspond to feasibility constraints associated with this adversarial realization. While this structure results in a more complex MP, the inclusion of second-stage variables and the use of primal cuts provides stronger approximations and faster convergence in practice, thereby improving overall computational performance \cite{113}.

  In this work, we implement a CCG algorithm that iteratively refines the approximation of the worst-case uncertainty realization by solving an MP and a set of adversarial SPs. However, before applying the decomposition strategy, it is necessary to reduce the original tri-level HARSO formulation to an equivalent two-level structure. This reformulation is achieved by using strong duality to replace the inner minimization, corresponding to the operational recourse problem, with its dual. In classical robust optimization settings, this transformation is valid when the inner problem is convex and satisfies regularity conditions, allowing the worst-case (max) and recourse (min) layers to be merged into a single maximization problem over the dual space. However, in the present case, the recourse problem includes a binary variable \( w_{j,t,y} \), which models the activation of charging and discharging operations in the battery system. To overcome this issue, we retain explicitly the binary variable \( w_{j,t,y} \) in the MP to be treated as a fixed parameter in the corresponding SP iteration. Once fixed, the SP becomes a convex linear program, allowing the application of strong duality and the construction of valid primal cuts that link the worst-case uncertainty realization to the operational cost.

  The SP is then formulated as the following primal SP:
\begin{equation}\label{Primal_SP_Model_Compact}
\begin{aligned}
\min\quad z =\; 
&\sum_{t\in \mathcal{T}}\sum_{y\in \mathcal{Y}}\sum_{s\in \mathcal{S}} \rho_s \Bigl (
    \lambda^{bg}_{t} \,p^{bg}_{t,y,s}
    -\lambda^{sg}_{t} \,p^{sg}_{t,y,s}
    +pg_{t,y,s} \,C_{\mathrm{op}}^{PV}\Bigr) \\
&+ \sum_{j\in \mathcal{J}}\sum_{t\in \mathcal{T}}\sum_{y\in \mathcal{Y}}\sum_{s\in \mathcal{S}}\,\rho_s (ds_{j,t,y,s}  \,C_{\mathrm{op}}^{BT} )
\end{aligned}
\end{equation}
\begin{subequations}\label{eq:sp constraints}
\begin{align}
\text{s.t.} \quad 
&pg_{t,y,s} - \widetilde{PL}_{t,y} - p^{sg}_{t,y,s} + p^{bg}_{t,y,s}
   + \sum_{j\in \mathcal{J}}\bigl(ds_{j,t,y,s}-ch_{j,t,y,s}\bigr) = 0 
&& \forall\,t,\,y,\,s &:\; \mathbf{a}_{t,y,s} \in \mathbb{R} \\
&pg_{t,y,s} \le \hat{\gamma}^{pv} \cdot \overline{PG}_{t,y,s} 
&& \forall\,t,\,y,\,s &:\; \mathbf{b}_{t,y,s} \leq 0 \\
&\text{soc}_{j,t,y,s} =
\begin{cases}
SOC^{i} \cdot \hat{\gamma}^{bt}_{j} + \varphi \cdot ch_{j,t,y,s} - \dfrac{1}{\varphi} \cdot ds_{j,t,y,s}, & t = 1 \\
SOC^{l} \cdot \hat{\gamma}^{bt}_{j}, & t = 24 \\
\text{soc}_{j,t-1,y} + \varphi \cdot ch_{j,t,y,s} - \dfrac{1}{\varphi} \cdot ds_{j,t,y,s}, & \text{otherwise}
\end{cases}
&& \forall\,j,\,t,\,y,\,s &:\; \mathbf{c}_{j,t,y,s} \in \mathbb{R} \\
&DG_{j,y} \cdot \hat{\gamma}^{bt}_{j} \cdot \underline{SOC} \le \text{soc}_{j,t,y,s}
&& \forall\,j,\,t,\,y,\,s &:\; \mathbf{d}^{-}_{j,t,y,s} \geq 0 \\
&\text{soc}_{j,t,y,s} \le DG_{j,y} \cdot \hat{\gamma}^{bt}_{j} \cdot \overline{SOC}
&& \forall\,j,\,t,\,y,\,s &:\; \mathbf{d}^{+}_{j,t,y,s} \leq 0 \\
&ch_{j,t,y,s} \le PB \cdot \hat{w}_{j,t,y,s}
&& \forall\,j,\,t,\,y,\,s &:\; \mathbf{f}_{j,t,y,s} \leq 0 \\
&ds_{j,t,y,s} \le PB \cdot (1 - \hat{w}_{j,t,y,s}) - PB \cdot (1 - \hat{\nu}^{bt}_j)
&& \forall\,j,\,t,\,y,\,s &:\; \mathbf{g}_{j,t,y,s} \leq 0
\end{align}
\end{subequations}

  The formulation shown above corresponds to the primal SP, which excludes first-stage constraints handled in the MP. All variables denoted with a hat symbol (e.g., \( \hat{\gamma}^{pv} \), \( \hat{w}_{j,t,y,s} \), \( \hat{\nu}^{bt}_j \)) are parameters passed from the MP, representing either investment decisions or fixed binary operational variables from the current iteration. The variable \( \widetilde{PL}_{t,y} \) denotes the worst-case realization of demand, resulting from the robust uncertainty set. All second-stage variables are scenario-dependent and indexed by \( s \in \mathcal{S} \), except for \( \widetilde{PL}_{t,y} \), which is fixed by the optimization of the adversarial SP. Dual variables are indicated to the right of each constraint, together with their corresponding sign restrictions, and will be used in the dual formulation to construct valid primal cuts. Thus, following the above dual notation, the dual SP is formulated as follows:
\begin{equation}\label{Dual_SP_Model}
\begin{aligned}
\max \quad z =\; 
&\sum_{t\in \mathcal{T}}\sum_{y\in \mathcal{Y}}\sum_{s\in \mathcal{S}}\left( 
    \widetilde{PL}_{t,y} \cdot \mathbf{a}_{t,y,s} 
    + \hat{\gamma}^{pv} \cdot \overline{PG}_{t,y,s} \cdot \mathbf{b}_{t,y,s} \right) \\
&+ \sum_{j\in \mathcal{J}}\sum_{y\in \mathcal{Y}}\sum_{s\in \mathcal{S}} \left(
    SOC^{i} \cdot \hat{\gamma}^{bt}_j \cdot \mathbf{c}_{j,1,y,s} 
  + SOC^{l} \cdot \hat{\gamma}^{bt}_j \cdot \mathbf{c}_{j,24,y,s} \right) \\
&+ \sum_{j\in \mathcal{J}}\sum_{t\in \mathcal{T}}\sum_{y\in \mathcal{Y}}\sum_{s\in \mathcal{S}} \Bigl(
    DG_{j,y} \cdot \hat{\gamma}^{bt}_j \cdot \underline{SOC} \cdot \mathbf{d}^{-}_{j,t,y,s} 
   \\ &+ DG_{j,y} \cdot \hat{\gamma}^{bt}_j \cdot \overline{SOC} \cdot \mathbf{d}^{+}_{j,t,y,s}
+ \mathbf{f}_{j,t,y,s} \cdot PB \cdot \hat{w}_{j,t,y,s} \\
&+ \mathbf{g}_{j,t,y,s} \cdot \left(PB \cdot (1 - \hat{w}_{j,t,y,s}) - PB \cdot (1 - \hat{\nu}^{bt}_j)\right) 
\Bigr)\\
\end{aligned}
\end{equation}
\begin{subequations}\label{eq:SP constraint 1}
\begin{align}
\text{s.t.} \quad
&\mathbf{a}_{t,y,s} + \mathbf{b}_{t,y,s} \leq \rho_s \cdot C_{\mathrm{op}}^{PV}
\quad  &&\forall\,t,\,y,\,s  \\
& \mathbf{a}_{t,y,s} \leq \rho_s \cdot \lambda^{bg}_{t}
\quad &&\forall\,t,\,y,\,s \\
& -\mathbf{a}_{t,y,s} \leq - \rho_s \cdot \lambda^{sg}_{t}
\quad &&\forall\,t,\,y,\,s \\
& -\mathbf{a}_{t,y,s} + \mathbf{f}_{j,t,y,s} \leq 0
\quad &&\forall\,j,\,y,\,s,\; t=24  \\
&-\mathbf{a}_{t,y,s} - \varphi \cdot \mathbf{c}_{j,t,y,s} + \mathbf{f}_{j,t,y,s} \leq 0
\quad &&\forall\,j,\,y,\,s,\; t<24 \\
& \mathbf{a}_{t,y,s} + \mathbf{g}_{j,t,y,s} \leq \rho_s \cdot C_{\mathrm{op}}^{BT}
\quad &&\forall\,j,\,y,\,s,\; t=24 \\
& \mathbf{a}_{t,y,s} + \frac{1}{\varphi} \cdot \mathbf{c}_{j,t,y,s} + \mathbf{g}_{j,t,y,s} \leq \rho_s \cdot C_{\mathrm{op}}^{BT}
\quad &&\forall\,j,\,y,\,s,\; t<24 \\
& \mathbf{d}^{-}_{j,t,y,s} + \mathbf{d}^{+}_{j,t,y,s} + \mathbf{c}_{j,t,y,s} \leq 0
\quad &&\forall\,j,\,y,\,s,\; t=24 \\
& \mathbf{d}^{-}_{j,t,y,s} + \mathbf{d}^{+}_{j,t,y,s} + \mathbf{c}_{j,t,y,s} - \mathbf{c}_{j,t+1,y,s} \leq 0
\quad &&\forall\,j,\,y,\,s,\; t<24
\end{align}
\end{subequations}

  We introduce binary variables \( V^{+}_{t,y}, V^{-}_{t,y} \) to integrate the BoU for the electricity demand into the dual formulation, since the optimal solution of the inner maximization problem lies at an extreme point of the uncertainty set. The uncertain parameter \( \widetilde{PL}_{t,y} \) is then redefined as an affine function of its respective nominal value and deviation bounds, constrained by the following expressions:
\begin{subequations}\label{eq:uncertainty_bou}
\begin{align}
&\widetilde{PL}_{t,y} = PL^{\mathrm{nom}}_{t,y} + \Delta PL_{t,y} \cdot V^{+}_{t,y} - \Delta PL_{t,y} \cdot V^{-}_{t,y} \label{eq:bou_pl}\\
&\sum_{t\in \mathcal{T}} (V^{+}_{t,y} + V^{-}_{t,y}) \leq \Gamma \qquad \forall y \in \mathcal{Y} \label{eq:bou_budget_pl}\\
&V^{+}_{t,y} + V^{-}_{t,y} \leq 1 \qquad \forall t \in \mathcal{T},\; \forall y \in \mathcal{Y} \label{eq:bou_binary_bound_pl}
\end{align}
\end{subequations}

Here, \( \Delta PL_{t,y} \) represents the maximum deviation from the nominal value. The parameter \( \Gamma \) controls the conservativeness of the uncertainty set, limiting how many time periods may simultaneously experience deviations from their nominal levels. However, when substituting expression~\eqref{eq:bou_pl} into the objective function of the dual SP~\eqref{Dual_SP_Model}, a bilinear term of the form \( V^{+}_{t,y} \cdot \mathbf{a}_{t,y,s} \) (and similarly \( V^{-}_{t,y} \cdot \mathbf{a}_{t,y,s} \)) appears. To handle this nonlinearity, a standard Big-M linearization approach is adopted \cite{mccormick1976computability}. This technique reformulates the product between a binary and a continuous variable, introducing auxiliary continuous variables ($\mathbf{d}^{-}, \mathbf{d}^{+}$) and additional constraints, ensuring the dual feasibility region remains linear and tractable. Therefore, the final dual representation of the SP using BoU consists of: (i) the dual problem of the inner operational SP~\eqref{Dual_SP_Model}, (ii) the Big-M constraints used to linearize the bilinear terms, and (iii) the BoU set constraints given in~\eqref{eq:bou_budget_pl} and~\eqref{eq:bou_binary_bound_pl}.

  The MP focuses on investment decisions while accounting for the worst-case electricity consumption and expected operational cost through the auxiliary variable \( \theta \). Thus, the MP is formulated as follows:
\setlength{\abovedisplayskip}{3pt}
\setlength{\belowdisplayskip}{3pt}
\begin{align}\label{FO_MP_model}
\min\quad z =\; & C^{PV}\, \gamma^{pv}
  \;+\;\sum_{j\in \mathcal{J}} C_{j}^{BT}\, \gamma^{bt}_{j} + \theta
\end{align}
\begin{subequations}\label{MP_Det_Const_model}
\begin{align}
\text{s.t.} \enspace   &\gamma^{pv}\;\le\;\overline{\Gamma^{PV}} \label{MP_PV_Cap_Const}\\
&w_{j,t,y}\;\le\;\nu^{bt}_{j}
   && \forall\,j\in \mathcal{J},t\in \mathcal{T},y\in \mathcal{Y}\label{MP_Binary_w_Const}\\
&\gamma^{bt}_{j}\;\le\;\nu^{bt}_{j}\;\cdot\;\overline{\Gamma^{BT}}
   && \forall\,j\in \mathcal{J\label{MP_BESS_Cap_Const}}\\
&\sum_{j\in \mathcal{J}}\nu^{bt}_{j}\;\le\;1 \label{MP_BESS_Tech_Const}\\
&\theta \geq {C_O^{k}}^\top y^k, & \forall & k = 1, \dots, r \label{Primal_cut_mp}\\
& h(x, y^k, u^k) = 0,               & \forall & k = 1, \dots, r \label{MP_Eq_const}\\
& g(x, y^k, u^k) \leq 0,            & \forall & k = 1, \dots, r \label{MP_Des_Const} 
\end{align}
\end{subequations}

Note that the term \( \theta \) represents the contribution of the second-stage operational cost under the worst-case demand realization. Constraints~\eqref{MP_PV_Cap_Const}--\eqref{MP_BESS_Tech_Const} define the first-stage decisions, including the binary variables \( w_{j,t,y} \). Constraint~\eqref{Primal_cut_mp} corresponds to the compact form of the primal cut, while Constraints~\eqref{MP_Eq_const}--\eqref{MP_Des_Const} capture the operational constraints that are included in the MP during iteration \( k \), based on the worst-case scenario \( k \) identified by the SP.

\begin{algorithm}[H]
\caption{CCG Algorithm for the HARSO problem}\label{alg:ccg_harso}
\begin{algorithmic}[1]
\State Initialize \( LB \gets -\infty \), \( UB \gets \infty \), tolerance \( \varepsilon \), iteration counter \( k \gets 0 \)
 \While{$\text{GAP} > \varepsilon$}
    \State \( k \gets k + 1 \)
    \State Solve the MP and obtain the current solution: \( \gamma^{pv}, \gamma^{bt}_j, \nu^{bt}_j \), and \( w_{j,t,s,y} \quad \forall j,t,s,y \)
    \State LB update: \( LB \gets z^{MP} \)
    \State Update the SP parameters: \( \hat{\gamma}^{pv} \gets \gamma^{pv} \),\quad \( \hat{\gamma}^{bt}_j \gets \gamma^{bt}_j,\quad \; \hat{\nu}^{bt}_j\gets \nu^{bt}_j \),\quad \ \(\hat{w}_{j,t,s,y}  \gets w_{j,t,s,y} \quad \forall j,t,s,y \)
    \State Solve the Dual SP and obtain optimal value \( z^{SP} \) and the worst-case demand \( \widetilde{PL}_{t,y} \)
    \State Update parameters in MP: \( \widehat{PL}_{t,y} \gets \widetilde{PL}_{t,y} \quad \forall t,y \)
    \State Add primal cut and constraints to MP:
\Statex \quad \( 
    \theta \geq \sum_{t,y,s} \rho_s \left( \lambda^{bg}_t \cdot p^{bg}_{t,y,s} 
    - \lambda^{sg}_t \cdot p^{sg}_{t,y,s} 
    + C_{\mathrm{op}}^{PV} \cdot pg_{t,y,s} \right) 
    + \sum_{j,t,y,s} \rho_s \cdot C_{\mathrm{op}}^{BT} \cdot ds_{j,t,y,s} 
\)
\Statex \quad \( 
    h(x, y^k, u^k) = 0, \enspace g(x, y^k, u^k) \leq 0
\)
    \State Update upper bound: \( UB \gets \min(UB, LB - \theta + z^{SP}) \)
    \State Compute relative gap: \( GAP \gets \left|1 - \frac{LB}{UB} \right| \)
\EndWhile
\end{algorithmic}
\end{algorithm}

Algorithm~\ref{alg:ccg_harso} outlines the CCG algorithm used to solve the proposed HARSO formulation. In Step~1, the algorithm is initialized by setting the lower and upper bounds (\( LB \leftarrow -\infty \), \( UB \leftarrow +\infty \)), defining a convergence tolerance \( \varepsilon \), and initializing the iteration counter \( k = 0 \). At each iteration, Step~4 solves the MP to determine the current investment decisions, namely the PV and BESS capacities \( \gamma^{pv}, \gamma^{bt}_j \), the battery technology selection \( \nu^{bt}_j \), and the operational schedule \( w_{j,t,s,y} \), which represent the system's long-term design and daily operation under the current information set. In Step~5, the lower bound is updated as \( LB \gets z^{MP} \), where \( z^{MP} \) denotes the objective function value of the MP. Next, in Steps~6--7, the optimal first-stage decisions are fixed and passed to the SP as parameters: \( \hat{\gamma}^{pv}, \hat{\gamma}^{bt}_j, \hat{\nu}^{bt}_j, \hat{w}_{j,t,s,y} \). The dual SP is then solved to obtain the worst-case realization of the demand \( \widetilde{PL}_{t,y} \) and its associated operational cost \( z^{SP} \). In Step~8, the worst-case demand is integrated into the MP, and the corresponding primal cut is added. Additionally, the operational constraints associated with this worst-case scenario \( h(x, y^k, u^k) = 0 \) and \( g(x, y^k, u^k) \leq 0 \) are appended to the MP. Finally, in Step~9, the upper bound is updated as \( UB \gets \min(UB, LB - 0 + z^{SP}) \), and in Step~10 the relative gap is computed using \( \text{GAP} \gets \left|1 - \frac{LB}{UB} \right| \). The algorithm iterates until the stopping criterion \( \text{GAP} < \varepsilon \) is satisfied.

\section{Case study and Computational results}\label{sec:case study}
This section presents the case study developed to evaluate the proposed planning model for PV and BESS systems under uncertainty. The study is conducted over a 10-year horizon, using historical irradiance and demand data to generate representative daily profiles for each year. The computational results assess the impact of different battery technologies, uncertainty levels, and system configurations on investment decisions and operational performance.

\subsection{Forecast of energy demand} \label{sec: energy demand}
The CAGR formula in (\ref{CAGR}) forecasts energy demand over the next 10 years. In this regard, we estimated the annual growth rate $r$ using data from \cite{17}, which provides projected energy demand values $D_t$ for 2025, 2030, and 2040, along with the corresponding base-year consumption $D_0$ and time intervals $t$. These projections account for the expected increase in energy demand due to future electrification trends in the country, including the installation of additional EV chargers in residential homes, the adoption of more efficient appliances, and other electrification-related developments. Using these values, $r$ was calculated by rearranging equation (\ref{CAGR}).
\begin{align}
    D_t = D_0 \cdot(1+r)^{t}
    \label{CAGR}
\end{align}

Next, to determine the base-year consumption $D_0$ for our projections, we analyzed electricity consumption data from multiple Spanish households over a one-year period. Some of these profiles corresponded to second homes, showing usage only on weekends or during specific weeks of the year. Others were considered non-representative due to extreme consumption patterns unusually high. This was detected since their average annual consumption presented very high values, meaning that they corresponded to commercial profiles. We identified these profiles, which were not representative of typical household consumption patterns, and we classified them using data clustering techniques. Further data preprocessing was carried out by imputing missing values with the commonly used KNN imputer method \cite{knn}, and after filtering and imputing the data, only few profiles remained for further treatment.
The data was then segmented into different periods, distinguishing between profiles with consumption of the same magnitude (constant) and those without. For profiles lacking constant consumption, a cyclic projection method was applied to reconstruct regular patterns, capturing multiple types of periodic seasonality (daily, weekly, and annual) as proposed by \cite{TAYLOR2010139}. 
From this processed dataset, a final representative profile $D_0$ was selected from a pool of 19 candidates, excluding homes equipped with photovoltaic panels, whose electricity bills do not accurately reflect total consumption.

Finally, $D_t$ projections of the obtained representative profile $D_0$ were made for the next 10 years ($t=1,2,...,10$) using the previously computed constant $r$ and the standard CAGR formula presented in (\ref{CAGR}). Hourly purchase prices were taken from \cite{price_data} and assumed constant over the following nine years.

\subsection{Estimated battery degradation} \label{Battery degradation}
Given the expected energy demand, we solved a deterministic optimization problem, referred to as the Self-Consumption model. This problem was used to determine the optimal operation of the battery within the PV-BESS system. In this model, uncertainties in energy generation and consumption are not considered for the sake of simplicity and because the primary objective is to approximate the battery’s operational behavior. Therefore, energy demand and production are the ones that were previously obtained with the CAGR and ARMA model. Additionally, the PV-BESS capacity is assumed based on the user's peak active power load. Based on this assumption, the following optimization problem is formulated and solved. The result is the hourly SOC level of the battery, obtained for every hour across the 10 years, given the expected PV generation and projected energy consumption.

We estimated the degradation of the PV-BESS batteries using the open-source BLAST software \cite{BLAST}, which provides a library of degradation and lifetime models for commercial lithium-ion batteries. Each battery model relies on publicly available lab-based aging data and considers parameters such as application, battery type, ambient temperature, and thermal management. The application is defined by the load profile, which in this study corresponds to the projected residential demand. Thus, the program can be applied to both FL and SL batteries, distinguished by their SOH, and incorporates experimental data that include calendar and accelerated aging. All models are based on lithium-ion technologies with different cathode chemistries, including Lithium iron phosphate (LFP), Lithium manganese oxide (LMO), Lithium nickel manganese cobalt oxide (NMC), and anode materials such as Graphite (Gr), Lithium titanate (LTO). The battery types considered in this study are summarized in Table \ref{tab:BLAST modelos de bateria}.

\begin{table}[htbp]
\centering
\begin{tabular}{|c|c|c|c|c|p{8cm}|}
\hline
\textbf{Type (Cathode/Anode)} & \textbf{Battery life stage} & \textbf{Capacity [Ah]} &\textbf{Lab data} &\textbf{Internal reference}\\
\hline
LFP/Gr & FL & 250 & \cite{lfp/gr} & Battery 1\\
LMO/Gr & SL & 66 & \cite{6}& Battery 2\\
NMC/Gr & FL & 75 &\cite{lfp/gr}& Battery 3\\
NMC/LTO & FL & 10 &\cite{nmc/lto}& Battery 4\\

\hline
\end{tabular}
\caption{Battery types included in BLAST \cite{BLAST} and considered in this study.}
\label{tab:BLAST modelos de bateria}
\end{table}

Ambient temperature is another input parameter of the program, as it significantly affects battery degradation. In this study, a one-year ambient temperature profile was established and assumed constant across the 10-year horizon, using meteorological data from \cite{climaSpain}. The specific estimation technique varies slightly across battery models due to the sensitivity of different cell types \cite{lfp/gr}, but all models account for key stressors that drive degradation. These include anode potential versus lithium (U\textsubscript{a}), which has been shown to influence calendar ageing \cite{105}, as well as temperature, SOC, depth of discharge (DoD), charge/discharge rate (C-rate), elapsed time, and equivalent full cycles (EFC). Further details on implementing degradation models for each chemistry are available in \cite{BLAST}.

\subsection{Data processing}
\begin{figure}[htbp]
    \centering
    \begin{subfigure}[b]{0.45\textwidth}
        \centering
        \includegraphics[width=\textwidth]{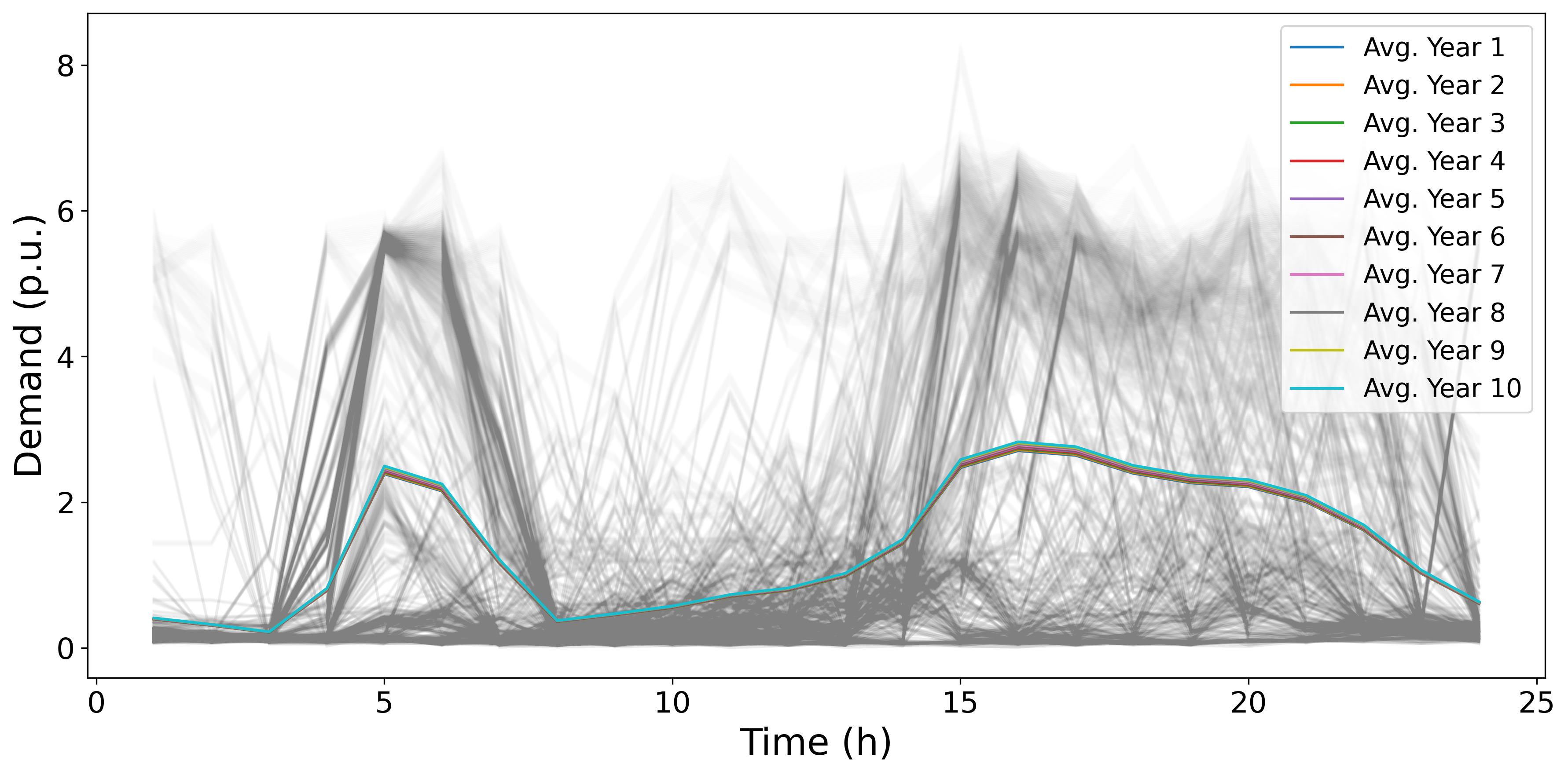}
        \caption{Estimated demand over 10 years.}
        \label{fig:PV generation over 10 years}
    \end{subfigure}
    \begin{subfigure}[b]{0.45\textwidth}
        \centering
        \includegraphics[width=\textwidth]{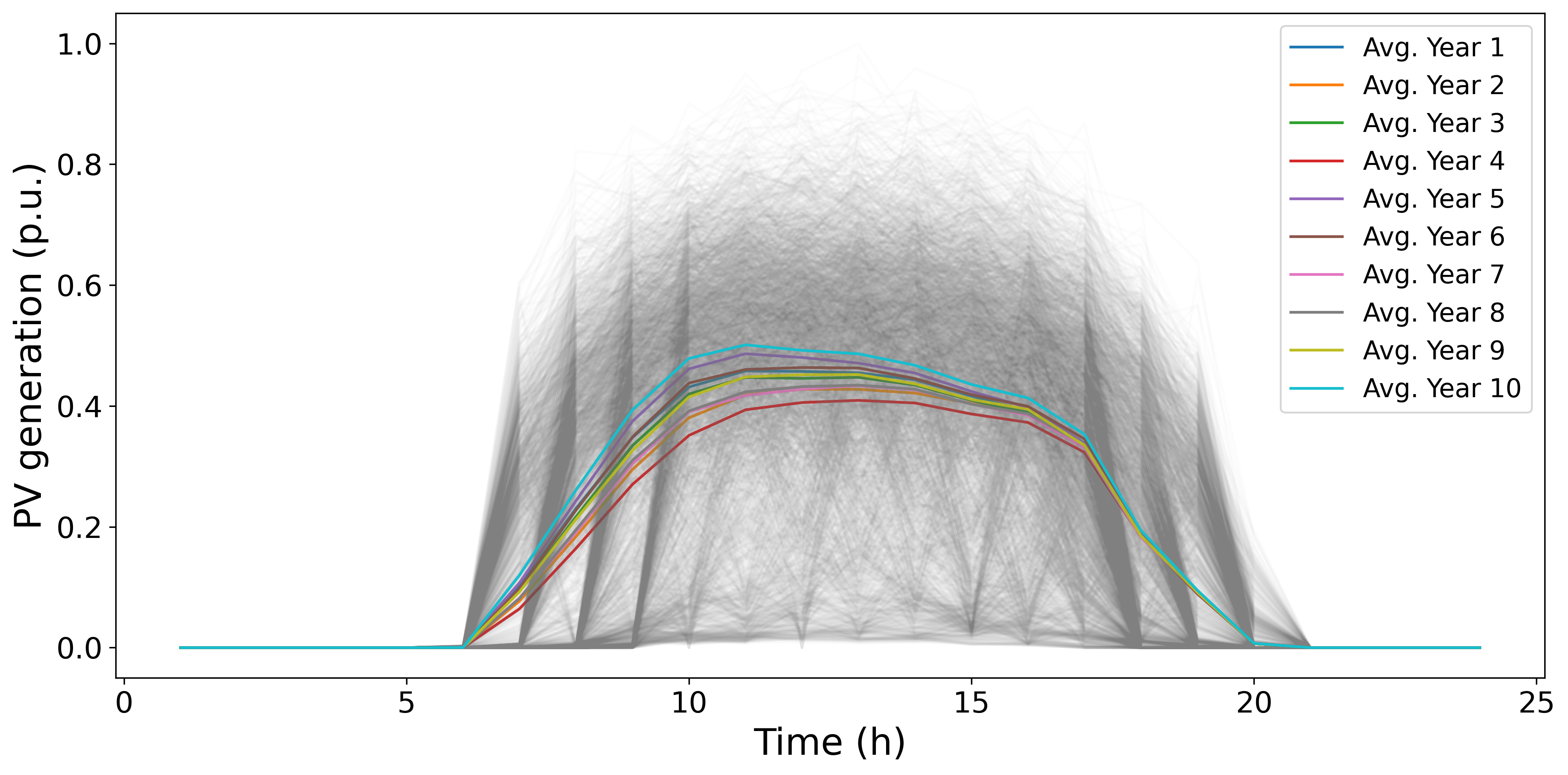}
        \caption{Estimated PV generation over 10 years.}
        \label{fig:estimated demand over 10 years}
    \end{subfigure}
    \caption{Estimated PV generation and energy demand over the next 10 years.}
    \label{fig:PV-and-demand 10 years}
\end{figure}

This study considers a representative rural household in Spain, with projected 10-year profiles of energy demand and PV generation shown in Figures \ref{fig:estimated demand over 10 years} and \ref{fig:PV generation over 10 years}, respectively.

Based on these projections, the corresponding battery SOC profiles are simulated across the full 10-year horizon using the deterministic simplified self-consumption model. Figure~\ref{fig:battery_operation} provides a closer look at the battery dynamics over a representative year, depicting the daily SOC trajectories. The figure highlights the average behavior (red line), the 5th and 95th percentiles (blue and green lines, respectively), and the full range of daily profiles (in gray), capturing the variability and typical operating patterns throughout the year. The shaded area represents the variability between the 5th and 95th percentiles, providing insight into the operational flexibility of the system. A consistent pattern is observed throughout the year: the battery typically charges during the afternoon, when electricity prices are lower, and discharges during periods of higher prices or lower PV generation, mostly in the evening. Notably, the model does not enforce a fixed final SOC at the end of each day, unlike typical approaches that use 24-hour representative scenarios. As a result, we observe that the average initial and final SOC values tend to stabilize between 20\% and 30\%, suggesting that this range could be considered efficient and realistic when defining boundary conditions in reduced-horizon scenarios. With this SOC profile, degradation is estimated for each of four batteries of study. In Figure \ref{fig:degradation_all_batteries}, a representation of the evolution of each battery SOH is shown. Section \ref{Battery degradation} presents the specifications of each battery material along with its corresponding internal reference. Battery 2 begins with a SOH of 70\% since it is a SL battery, while the remaining batteries, being FL, start at 100\% SOH. As it can be seen, degradation patterns vary significantly across chemistries, and the SL battery exhibits a faster decline in SOH compared to the FL batteries. This further highlights, once more, the relevance of incorporating multiple battery chemistries into the analysis, given their differing operational and degradation characteristics.

\begin{minipage}[t]{0.48\textwidth}
  \begin{figure}[H]
    \centering
    \includegraphics[width=\textwidth]{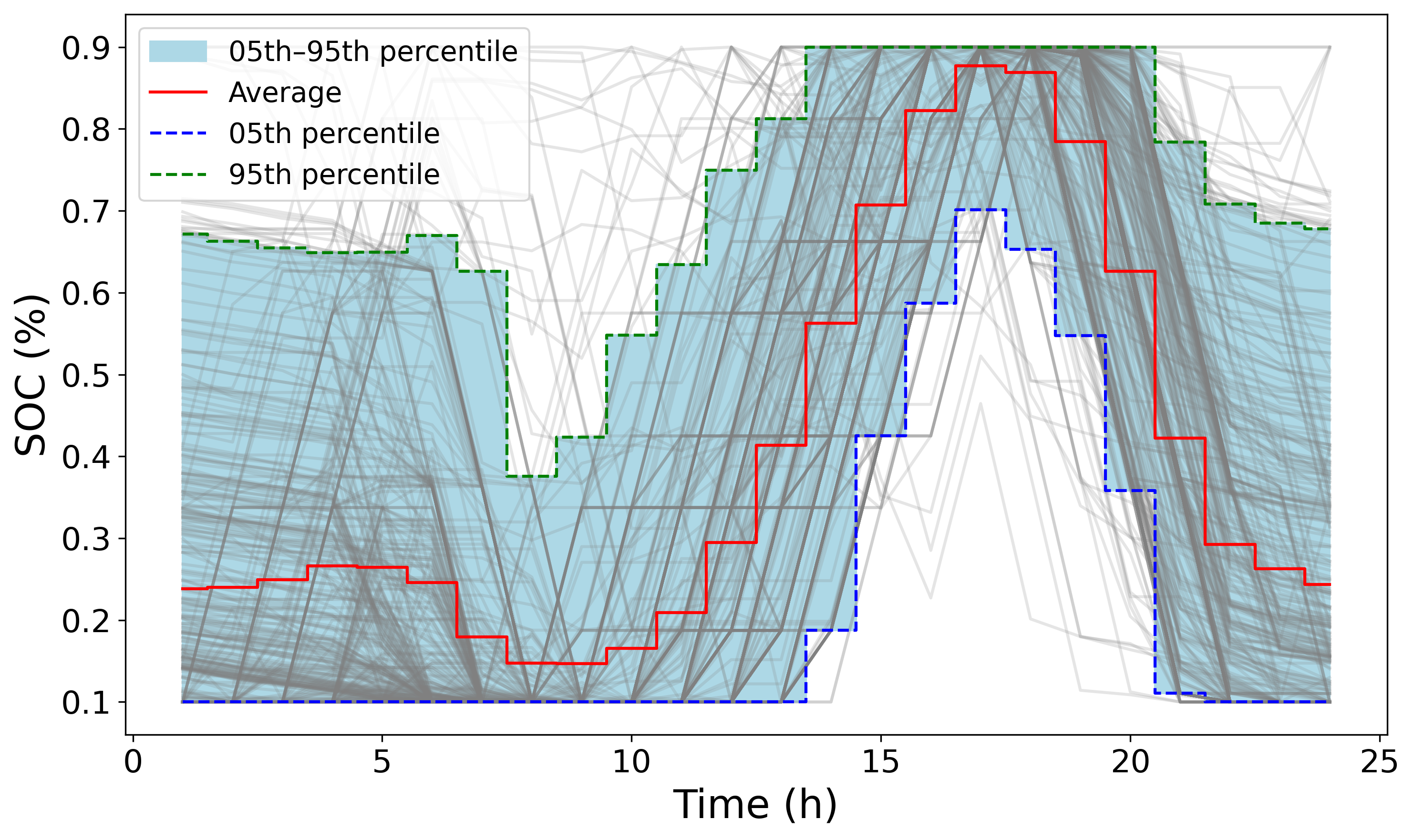}
    \caption{Estimated battery operation over time.}
    \label{fig:battery_operation}
  \end{figure}
\end{minipage}
\hfill
\begin{minipage}[t]{0.48\textwidth}
  \begin{figure}[H]
    \centering
    \includegraphics[width=\textwidth]{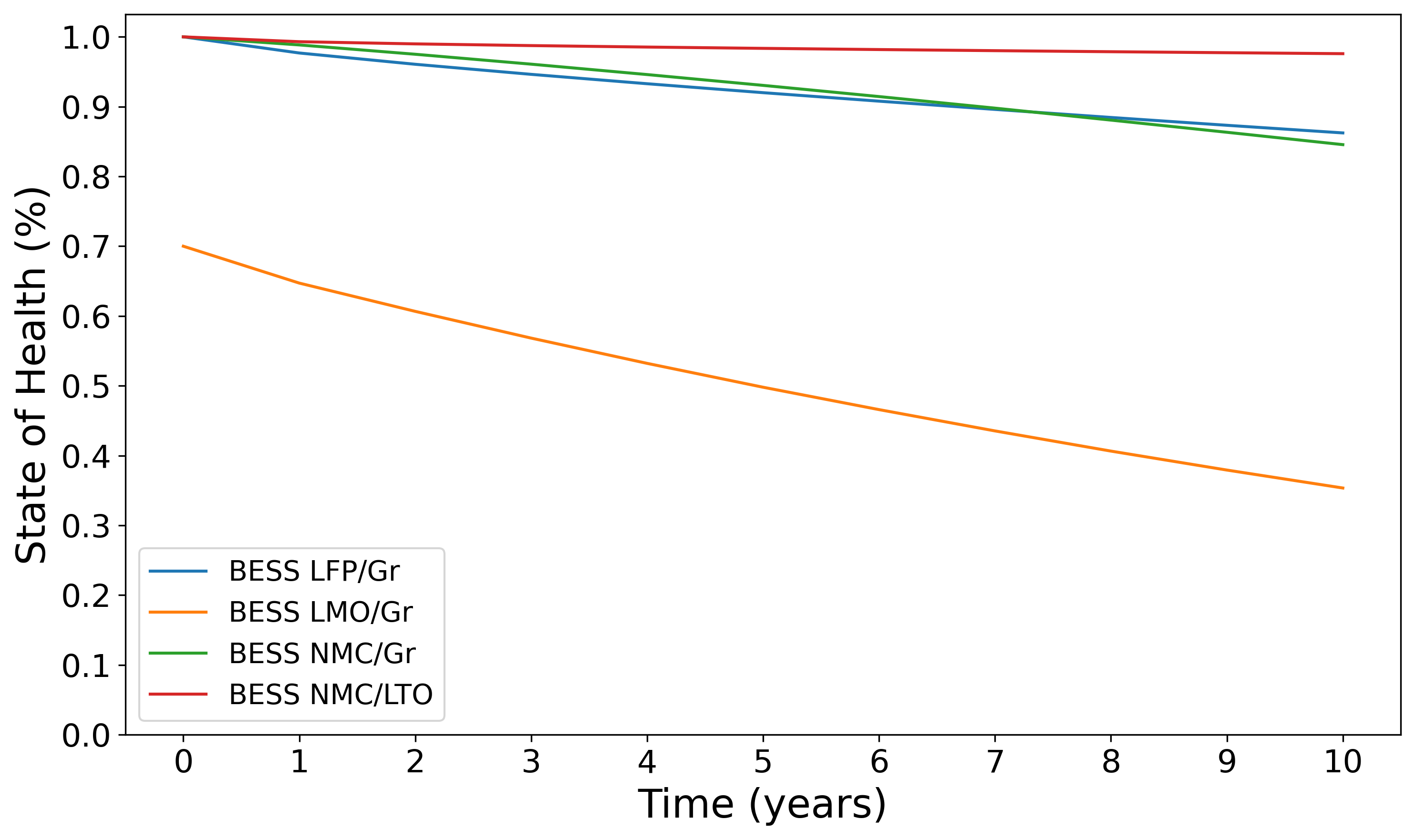}
    \caption{Degradation profiles over time for all batteries.}
    \label{fig:degradation_all_batteries}
  \end{figure}
\end{minipage}
 
\vspace{0.5cm}


Representative 24-hour profiles of generation and demand for each projected year have been used to incorporate uncertainty in the HARSO model. These representative profiles are derived from the 10-year hourly forecasts. For PV generation, four distinct scenarios are considered, potentially representing seasonal variability (e.g., spring, summer, autumn, and winter). Figure~\ref{fig:PV_generation} shows the hourly PV generation profiles for year 10 under these four scenarios, highlighting the differences in solar availability across the year. For energy demand, hourly values are used to compute the average, minimum, and maximum consumption levels observed in each projected year. These statistics are then used to define the uncertainty interval for demand. Figure~\ref{fig:power_demand_y10} displays the 24-hour demand profile for year 10, showing the average (red), minimum (blue), and maximum (green) curves, along with the full distribution of daily demand patterns in the background.

\begin{minipage}[t]{0.48\textwidth}
  \begin{figure}[H]
    \centering
    \includegraphics[width=\textwidth]{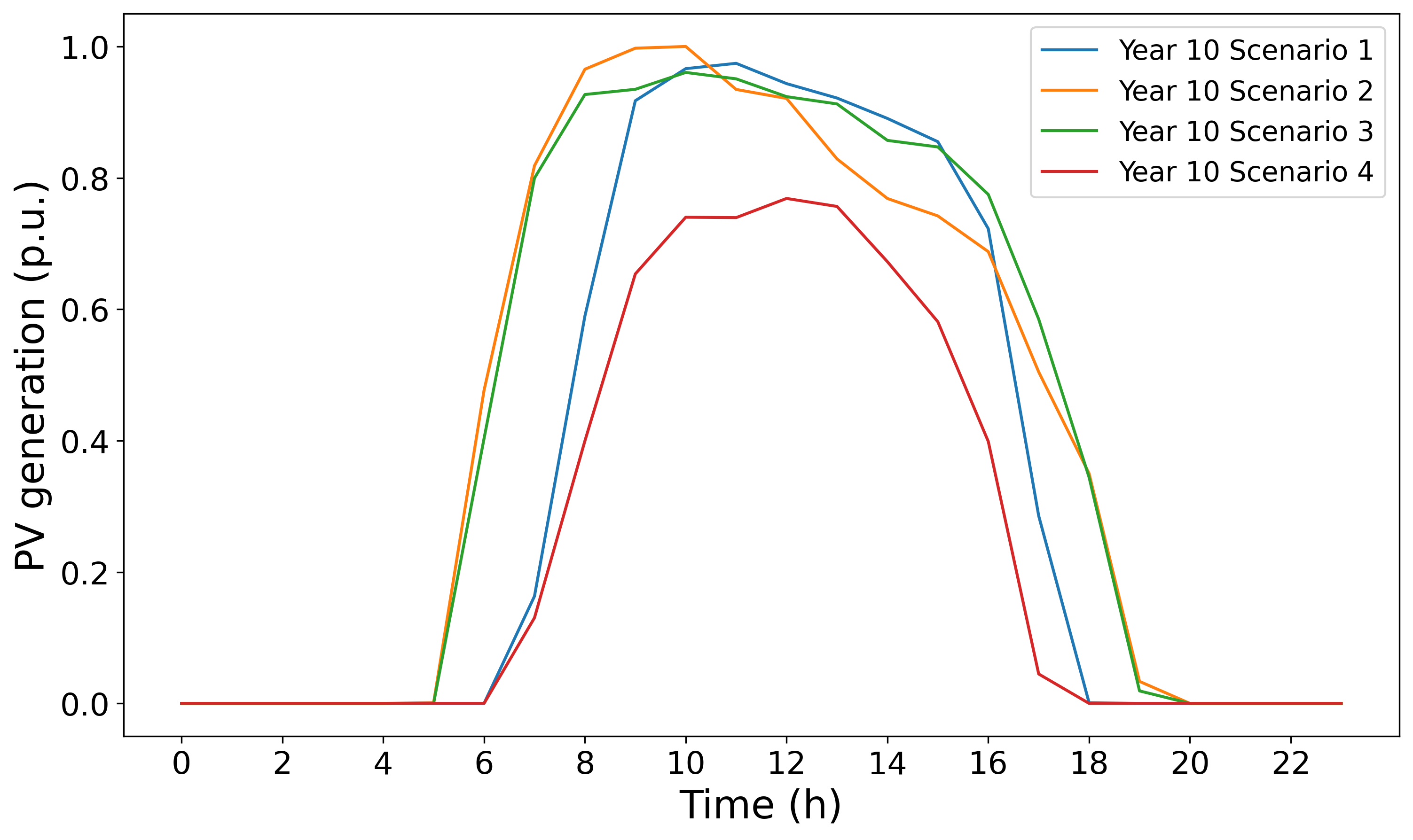}
    \caption{24-hour PV generation under 4 scenarios at year 10.}
    \label{fig:PV_generation}
  \end{figure}
\end{minipage}
\hfill
\begin{minipage}[t]{0.48\textwidth}
  \begin{figure}[H]
    \centering
    \includegraphics[width=\textwidth]{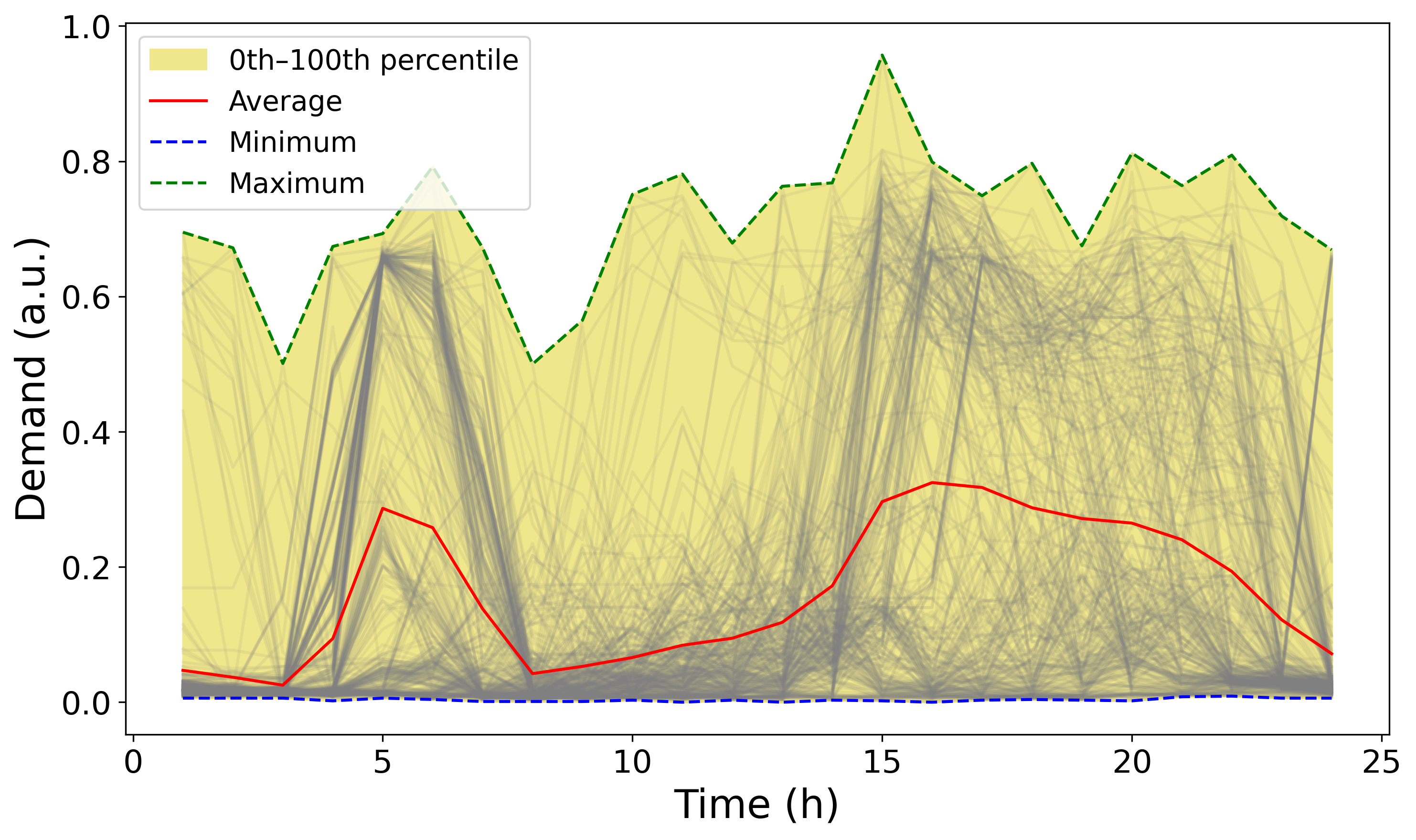}
    \caption{24-hour power demand intervals at year 10.}
    \label{fig:power_demand_y10}
  \end{figure}
\end{minipage}

\vspace{0.5\baselineskip}



\subsection{Sensitivity analysis on the BoU} \label{sensitivity_analysis_scenarios}
Battery investment costs were derived from current commercial sources, with Battery 2 being the least expensive, followed by Batteries 1 and 3, and Battery 4 being the most costly. Yet, when considering both capacity and performance in a simple cost analysis, Battery 1 shows the best economic performance from an operational perspective, followed by Batteries 4, 3, and 2. Thus, with an initial BoU of 5, selected initially as a balanced value that introduces moderate uncertainty without being overly conservative, the optimization results indicate that Battery 1 (LFP/Gr chemistry) is the most cost-effective option. Under this setting, the optimal investment corresponds to an installed PV capacity of 13.6 kWh and a BESS capacity of 6.10 kWh, resulting in a total cost of \EUR{61.95} per day.

\begin{figure}[htbp]
    \centering
    \includegraphics[width=0.6\linewidth]{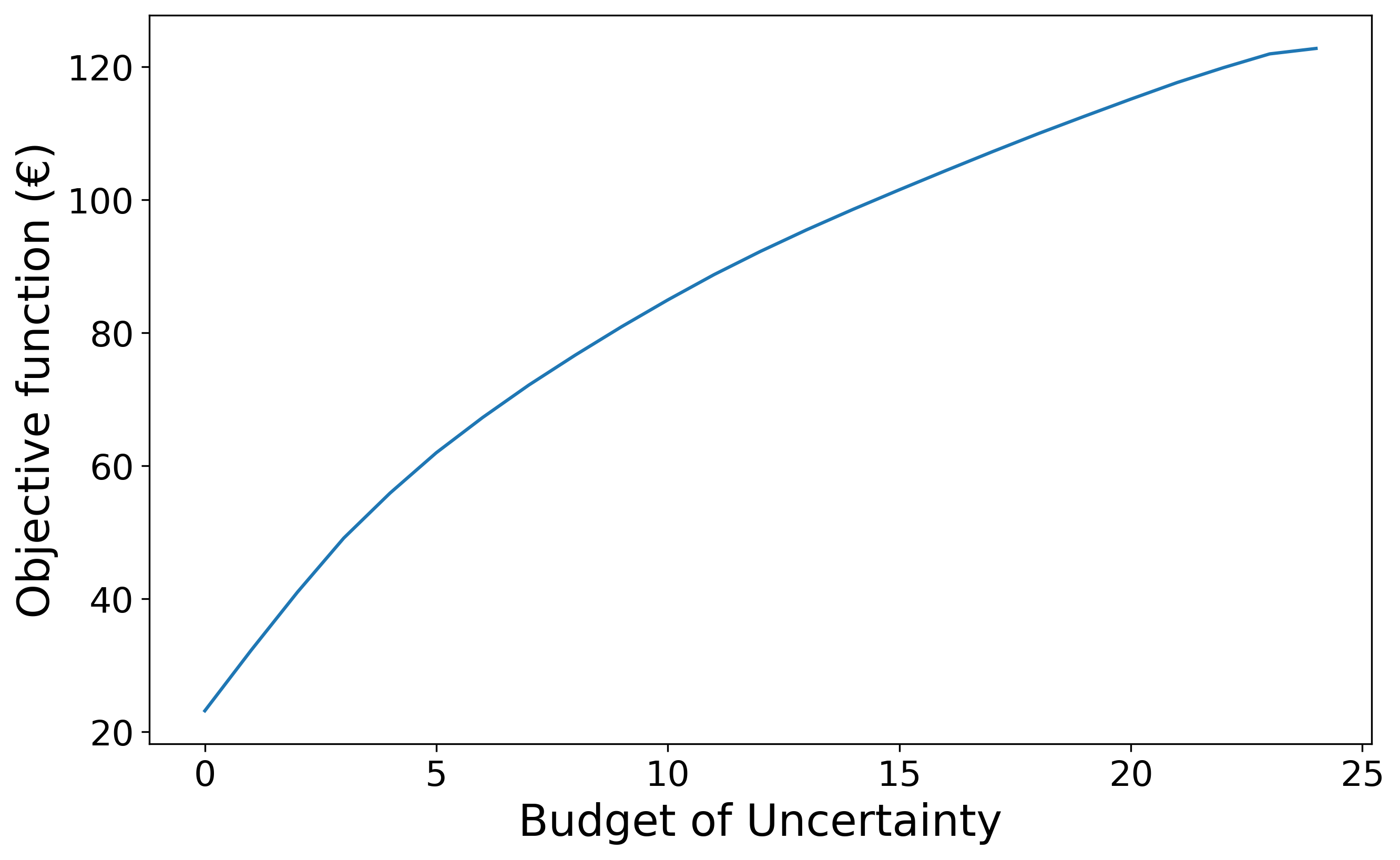}
    \caption{Objective function evolution for multiple BoU and 4 scenarios.}
    \label{fig:FunObj_24BoU_4Scenarios}
\end{figure}

Given that the BoU parameter regulates the degree of protection against worst-case demand deviations, ranging from 0 (fully deterministic) to 24 (fully robust, where all deviations are considered simultaneously), a sensitivity analysis was performed by varying the arbitrary selection of its initial BoU. Thus, Figure~\ref{fig:FunObj_24BoU_4Scenarios} presents the evolution of the minimized objective function as BoU increases, considering four PV generation scenarios. The curve exhibits a concave shape, meaning that the marginal increase in cost diminishes as BoU grows. This indicates that most of the cost impact occurs at lower BoU values, while higher levels of robustness lead to smaller additional investments. 

The impact of the number of scenarios on the objective function was also tested. Results showed that increasing the number of PV generation scenarios (e.g., from four to six) had minimal influence on the evolution of the objective function, yielding nearly identical curves. For this reason, only the case with four scenarios is presented in the figure, as additional plots would replicate the same trend.

\begin{figure}[htbp]
    \centering
    \includegraphics[width=0.65\linewidth]{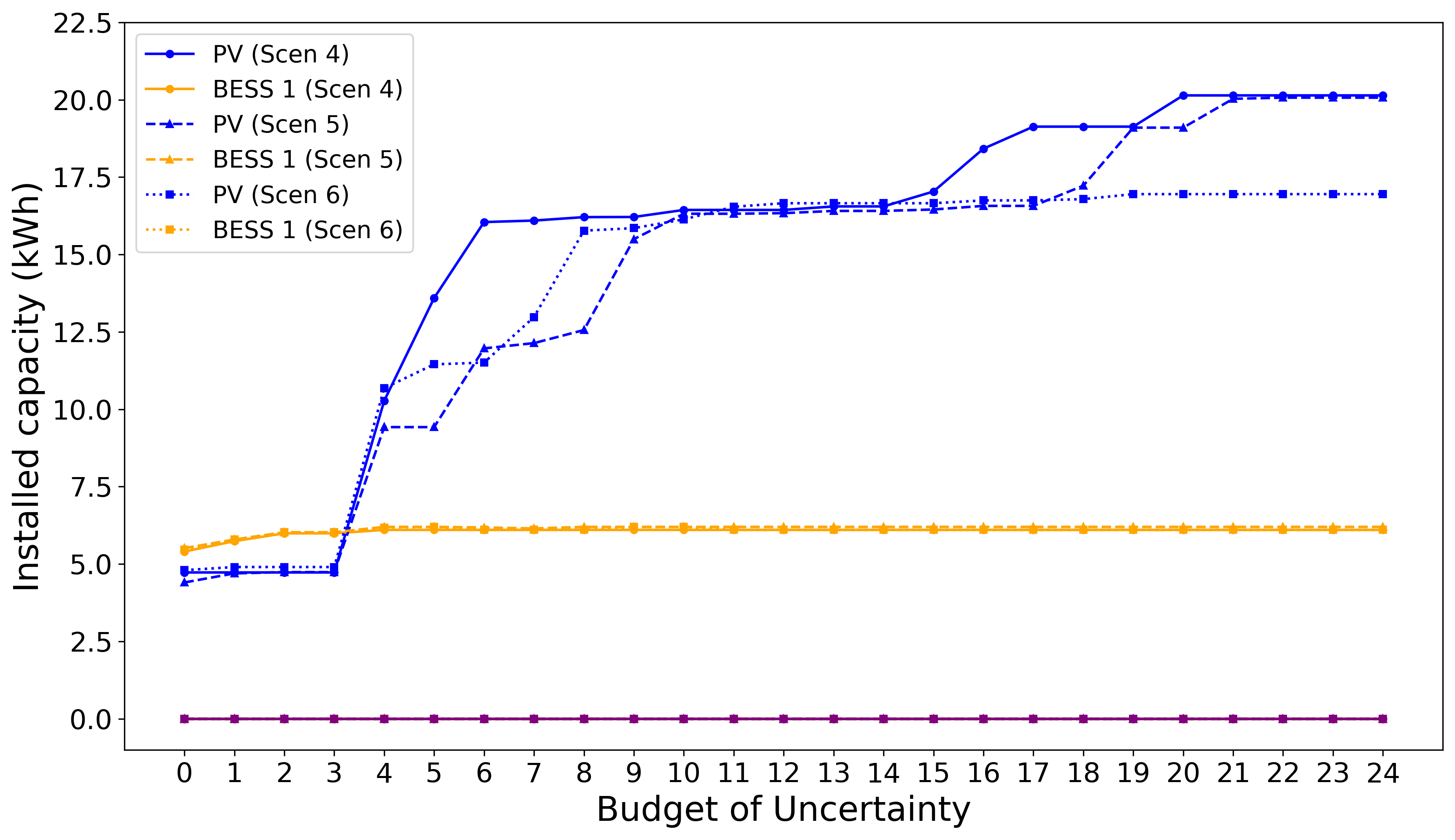}
        \caption{PV-BESS optimal installation under different ranges of BoU and 4, 5 and 6 scenarios.}
        \label{fig:PV_BESS_BoU_24_4scenarios_Comparison}
\end{figure}

Figure~\ref{fig:PV_BESS_BoU_24_4scenarios_Comparison} shows the evolution of the optimal PV and BESS installed capacities as the BoU increases in instances with four, five and six PV generation scenarios. As illustrated, both technologies converge to stable values, although they follow different trajectories. In this regard, for low BoU values (BoU $\leq$ 3), the installed PV capacity remains constant, and robustness is achieved primarily by installing more BESS (FL type) capacity than PV. However, as BoU increases,  results show that further increasing BESS capacity becomes suboptimal, and increasing the PV installation is more beneficial. Interestingly, there is a clear shift in how the model responds to uncertainty. With low values of BoU, more BESS is cost-effective. With greater values of BoU, more PV becomes preferable. This reflects a strategic shift in design philosophy, driven by the conservativeness of energy demand assumptions: from enhancing operational flexibility through storage, to mitigating uncertainty through generation overcapacity. Likewise, BESS type 1 consistently emerges as the most suitable option across all levels of uncertainty. This result holds when the number of PV scenarios is increased to five and six, reinforcing its dominance. Thus, FL BESS's favorable degradation profile and relatively low investment and operational costs explain its consistent selection.

Although the number of PV generation scenarios has a limited impact on the overall system configuration, some differences can be observed depending on the BoU range. In proportional terms, increasing the number of scenarios from four to six does not substantially alter the convergence value of the installed BESS capacity nor the BoU threshold at which PV becomes more cost-effective than storage. The overall design strategy remains largely consistent. However, a closer examination of Figure~\ref{fig:PV_BESS_BoU_24_4scenarios_Comparison} reveals that the number of PV generation scenarios can influence the installed capacity depending on the BoU range. Between BoU values of 4 and 9, increasing the number of scenarios from 4 to 6 leads to slightly lower installed PV capacities. This suggests that a richer scenario set under moderate uncertainty allows the model to distribute generation variability better and avoid overinvestment. However, for intermediate BoU values (between 10 and 16), all scenario sets yield nearly identical PV capacities, indicating a potential convergence zone where the design appears to stabilize across the considered scenario granularities. Interestingly, when BoU exceeds 18, the case with six scenarios diverges again from the others, consistently recommending lower PV capacity. This behavior may reflect a better anticipation of extreme seasonal variations in solar output, allowing the model to limit overdimensioning even under highly conservative settings.

\begin{figure}[htbp]
    \centering
    \includegraphics[width=0.6\linewidth]{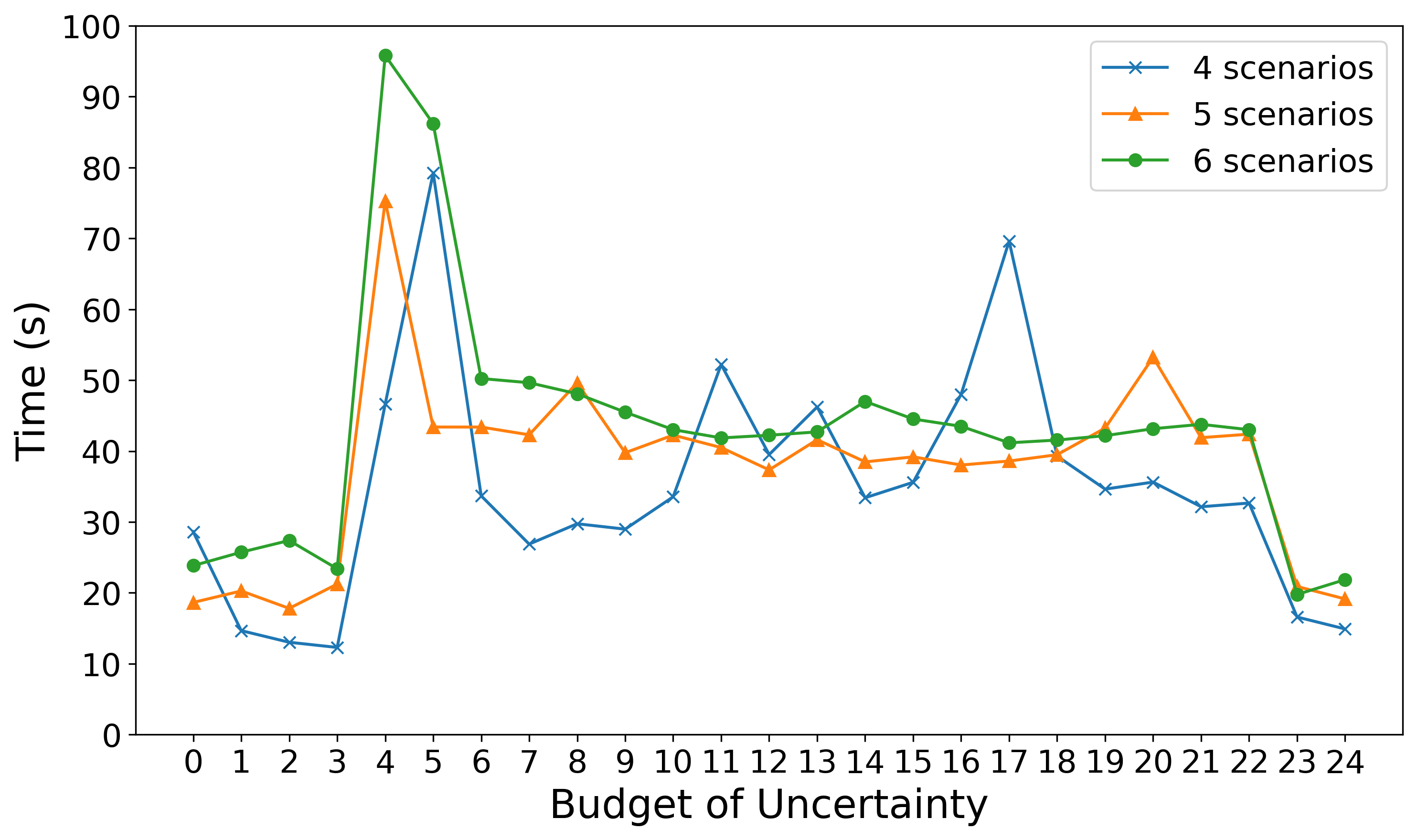}
    \caption{Required computational time for multiple BoU and scenarios.}
    \label{fig:time_BoU_scenarios}
\end{figure}

Figure~\ref{fig:time_BoU_scenarios} illustrates the computational time required to solve the HARSO model for different BoU values and scenario sets. As expected, models with a higher number of scenarios tend to require more time on average, due to the increased dimensionality and size of the underlying optimization problem. However, the relationship between BoU and computational time is not linear. A consistent computational peak is observed around BoU values of 4 and 5 across all scenario sets, as observed in Figure~\ref{fig:time_BoU_scenarios}. This suggests that intermediate levels of uncertainty introduce greater complexity for the solution algorithm. In this range, the model must explore a broader and less structured space of adversarial combinations, likely involving strategic shifts in technology deployment, such as transitioning from BESS-based to PV-based robustness. In contrast, for very low BoU values, the problem becomes nearly deterministic and thus easier to solve. Similarly, the system's response stabilizes for high BoU values (BoU $\geq$ 18), leading to faster convergence as decisions remain relatively unchanged across iterations.

\subsection{Trade-offs between uncertainty and cost in BoU selection}


As we have seen before, adding uncertainty considerations naturally leads to higher computational efforts. In a similar way, it can also lead to more expensive system configurations. In this section we analyze this trade-off for our case study.

Based on the analysis presented in Section~\ref{sensitivity_analysis_scenarios}, increasing the number of PV generation scenarios does not substantially alter the system’s optimal configuration. Although a larger scenario set may offer a more detailed representation of uncertainty, it also leads to a noticeable increase in computational time. Therefore, using four scenarios is assumed sufficient to balance accuracy and efficiency in the following analysis.

\begin{figure}[htbp]
    \centering
    \includegraphics[width=0.6\linewidth]{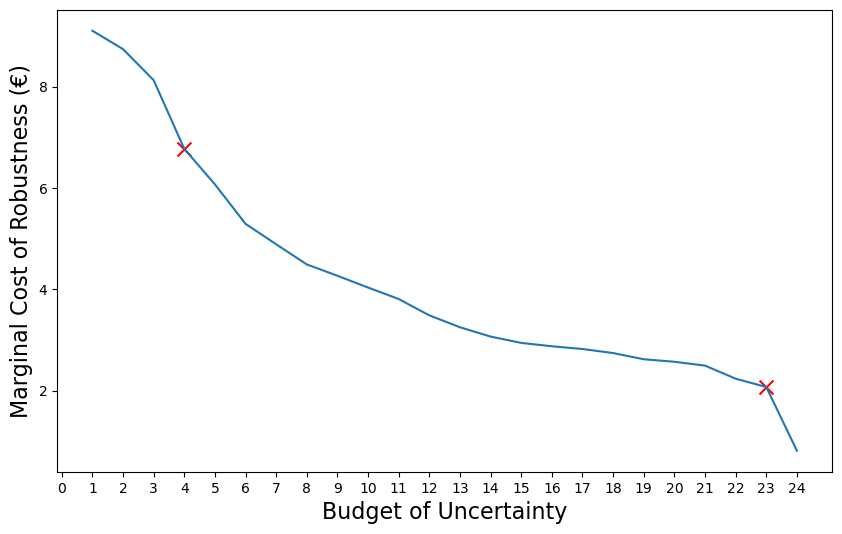}
    \caption{Marginal cost of robustness evolution.}
    \label{fig:MCR}
\end{figure}

On the other hand, the BoU depends on how much conservatism is desired. As the BoU increases, the level of conservatism increases, and the objective function also increases as seen in Figure~\ref{fig:FunObj_24BoU_4Scenarios}, i.e. more conservative solutions come at a higher cost. However, Figure~\ref{fig:MCR} shows that the marginal cost of robustness (MCR) is decreasingly concave, implying that the marginal cost of the added robustness decreases as the BoU grows. Therefore, adding more conservatism leads to higher costs, although the increase becomes smaller with higher conservatism. 
In Figure \ref{fig:MCR}, two inflection points are identified at BoU$=4$ and BoU$=23$. The second inflection point at BoU$=23$ is not very representative, as the trend beyond that point cannot be fully observed given the maximum value of BoU$=24$. From BoU$=4$ onward, the rate of decrease in marginal cost begins to slow. This indicates that for BoU$<4$, the level of conservatism must be chosen carefully, as the marginal cost of allowing one extra time period to be uncertain is very high. Past this point, adding an additional time period becomes progressively less costly, allowing the decision-maker to evaluate trade-offs between cost and robustness more effectively. 

\section{Conclusions}\label{sec:conclusions}
In this study, a HARSO model was proposed to determine the optimal PV-BESS design. We model the uncertainty in the energy generation using stochastic optimization, and the uncertainty in the energy demand using robust optimization.
We considered different battery types, each with its corresponding estimated degradation profile. 
We used a CCG algorithm to solve the HARSO model, which provides a more efficient resolution compared to commonly used approaches in the literature, such as Benders decomposition.

To demonstrate the model’s applicability, we solved a case study involving a household in Spain. Results show that the optimal installed capacities are not significantly influenced by the number of scenarios considered. However, they are strongly affected by the level of conservatism, as defined by the BoU. A critical BoU threshold is identified, beyond which the model shifts its optimal strategy: for lower BoU values, it is optimal to increase battery capacity to ensure robustness against worst-case energy demand; for higher BoU values, the model favors increasing the installed PV capacity instead. A sensitivity analysis reveals that selecting BoU values below 4 must be considered very carefully by the decision-maker, as additional conservatism leads to a substantial increase in system costs. In contrast, for BoU$\geq4$, the trade-off between cost and robustness is more stable, allowing the decision-maker to balance more easily the additional robustness against cost increases. Finally, the same battery type consistently emerges as the optimal choice, due to its low degradation rate and relatively low investment and operational costs.

Further development is anticipated in future research efforts, particularly regarding improvements in solution methods. First, we wish to explore the implementation of Benders decomposition as a branch-and-cut strategy, along with other decomposition techniques aimed at simplifying the model. Additionally, we want to extend the framework to larger-scale problems involving multiple households or entire energy communities. More extreme scenarios for energy generation could also be considered, and the current 24-hour time window could be expanded to capture longer-term system behavior.

\section*{Acknowledgment}
This work has been supported by ANID FONDECYT Iniciación 11240745 and Fondequip Mediano EQM200234.

\bibliographystyle{APA_Elsevier}
\bibliography{references}


\newpage

\begin{table}[htbp]   
\small
\begin{mdframed}
\nomenclature[01]{\textbf{Sets}}{}
\nomenclature[02]{$h \in \mathcal{H}$}{Set of time slots in a year.}
\nomenclature[03]{$t \in \mathcal{T}$}{Set of hours in a day.}
\nomenclature[04]{$s \in \mathcal{S}$}{Set of uncertainty scenarios.}
\nomenclature[05]{$y \in \mathcal{Y}$}{Set of years.}
\nomenclature[06]{$j \in \mathcal{J}$}{Set of battery types.}
\nomenclature[07]{\textbf{Parameters}}{}
\nomenclature[08]{$PL_t$}{Power load at time $t$ [\%].}
\nomenclature[09]{$\overline{PV}_t$}{Maximum PV active power available to inject at period $t$ [kW].}
\nomenclature[10]{$\overline{SOC}$}{Maximum state of charge of BESS [kWh].}
\nomenclature[11]{$\underline{SOC}$}{Minimum state of charge of BESS [kWh].}
\nomenclature[12]{$\varphi$}{charging or discharging BESS efficiency [\%].}
\nomenclature[13]{$PB$}{Maximum power rate for BESS charging or discharging [kW].}
\nomenclature[14]{$\Gamma^{PV}$}{PV capacity installed [kW]}
\nomenclature[15]{$\Gamma^{BT}$}{BESS capacity installed [kWh]}
\nomenclature[16]{$\varphi$}{Charging/discharging efficiency of BESS [\%].}
\nomenclature[17]{$M$}{Big–\(M\) constant for linearization.}
\nomenclature[18]{$\overline{PL}$}{Maximum normalized load factor.}
\nomenclature[19]{$C^{PV}$}{Investment cost of PV system [\$/kW].}
\nomenclature[20]{$C_{\mathrm{op}}^{PV}$}{Operational cost of PV per period [\$/kW].}
\nomenclature[21]{$C_{\mathrm{op}}^{BT}$}{Operational cost of battery per period [\$/kWh].}
\nomenclature[22]{$\overline{\Gamma^{PV}}$}{Upper bound on installable PV capacity [kW].}
\nomenclature[23]{$\overline{\Gamma^{BT}}$}{Upper bound on installable BESS capacity [kWh].}
\nomenclature[24]{$C^{BT}_{j}$}{Investment cost of BESS system of type \(j\) [\$/kW].}
\nomenclature[25]{$PL_{t,y}$}{Normalized load at time \(t\) in year \(y\) [\%].}
\nomenclature[26]{$\overline{PG}_{t,s,y}$}{Maximum PV generation at time \(t\), scenario \(s\), year \(y\) [kW].}
\nomenclature[27]{$\lambda_{t}^{bg}$}{Power price bought from the energy supplier at time \(t\) [\$/kW].}
\nomenclature[28]{$\lambda_{t}^{sg}$}{Power price sold to the energy supplier at time \(t\) [\$/kW].}
\nomenclature[29]{$DG_{j,y}$}{Degradation rate of battery type \(j\) in year \(y\) [\%].}
\nomenclature[30]{$\rho_{s}$}{Probability of scenario \(s\).}
\nomenclature[31]{\textbf{Variables}}{}
\nomenclature[32]{$p^{bg}_{t}$}{Power bought from the energy supplier at time $t$ [kW]}
\nomenclature[33]{$p^{sg}_{t}$}{Power sold to the energy supplier at time $t$ [kW]}
\nomenclature[34]{$p^{PV}_{t}$}{Power from the PV system at time $t$ [kW]}
\nomenclature[35]{$ds_{t}$}{Discharged power from the BESS at time $t$ [kW]}
\nomenclature[36]{$ch_{t}$}{Charged power to the BESS at time $t$ [kW]}
\nomenclature[37]{$soc_{t}$}{BESS state of charge at time $t$ [kWh]}
\nomenclature[38]{$w_{t}$}{1 if BESS is charging at time $t$. 0 otherwise.}
\nomenclature[39]{$pg_{t,s,y}$}{Power generated by the PV panels at time \(t\), scenario \(s\) and year \(y\) [kW].}
\nomenclature[40]{$p^{bg}_{t,s,y}$}{Power bought from the energy supplier at time \(t\), scenario \(s\) and year \(y\)[kW].}
\nomenclature[41]{$p^{sg}_{t,s,y}$}{Power sold to the energy supplier at time \(t\), scenario \(s\) and year \(y\) [kW].}
\nomenclature[42]{$soc_{j,t,s,y}$}{State of charge of battery type \(j\) at time \(t\), scenario \(s\) and year \(y\) [kWh].}
\nomenclature[43]{$ch_{j,t,s,y}$}{Charging power to battery type \(j\) at time \(t\), scenario \(s\) and year \(y\) [kW].}
\nomenclature[44]{$ds_{j,t,s,y}$}{Discharging power from battery type \(j\) at time \(t\), scenario \(s\) and year \(y\) [kW].}
\nomenclature[45]{$w_{j,t,s,y}$}{Binary variable; 1 if battery \(j\) is charging at time \(t\), scenario \(s\) and year \(y\), 0 otherwise.}
\nomenclature[46]{$\gamma^{pv}$}{Installed PV capacity [kW].}
\nomenclature[47]{$\gamma^{bt}_{j}$}{Installed capacity of battery type \(j\) [kWh].}
\nomenclature[48]{$\nu^{bt}_{j}$}{Binary variable; 1 if battery type \(j\) is installed, 0 otherwise.}
\printnomenclature
\end{mdframed}
\caption{Nomenclature for the deterministic Self-Consumption model.}\label{Nomenclature_table}
\end{table}

\end{document}